\newcommand\BB{\mathbb B}
\newcommand\CC{{\mathbb C}}

\newcommand\cD{{\cal D}}
\newcommand\cE{{\cal E}}
\newcommand\cF{{\cal F}}
\newcommand\cG{{\cal G}}

\newcommand\cK{{\cal K}}
\newcommand\cL{{\cal L}}
\newcommand\cM{{\cal M}}

\newcommand\cO{{\cal O}}
\newcommand\cP{{\cal P}}
\newcommand\cQ{{\cal Q}}

\newcommand\cS{{\cal S}}

\newcommand\cU{{\cal U}}
\newcommand\cV{{\cal V}}

\newcommand\cX{{\cal X}}
\newcommand\cY{{\cal Y}}

\newcommand\es{\emptyset}

\newcommand\HH{{\mathbb H}}
\newcommand\hra{\hookrightarrow}
\newcommand\lra{\longrightarrow}

\newcommand\ov{\overline}
\newcommand\PP{{\mathbb P}}
\newcommand\QQ{{\mathbb Q}}

\newcommand\wh{\widehat}
\newcommand\wt{\widetilde}
\newcommand\ZZ{{\mathbb Z}}

\newcommand{\mapor}[1]{{\stackrel{#1}{\longrightarrow}}}
\newcommand{\mapver}[1]{\Big\downarrow
\vcenter{\rlap{$\scriptstyle#1$}}}
\documentclass[10pt,a4paper]{article}
\usepackage{amsthm}
\usepackage{amsmath}
\usepackage{amssymb}
\newtheorem{thm}{Theorem}[section]
\newtheorem{clm}[thm]{Claim}

\newtheorem{crl}[thm]{Corollary}
\newtheorem{dfn}[thm]{Definition}
 
\newtheorem{lmm}[thm]{Lemma}
\newtheorem{prp}[thm]{Proposition}

\begin{document}
 \title{Dual double EPW-sextics and their periods}
 \author{Kieran G. O'Grady\\
Universit\`a di Roma \lq\lq La Sapienza\rq\rq}
\date{April 12 2006
\vskip 4mm \small{\it Dedicato a Fedya Bogomolov in occasione del $60^{\text{o}}$ compleanno}}
 \maketitle
 \section{Introduction}\label{prologo}
 \setcounter{equation}{0}
Eisenbud, Popescu and Walter have constructed certain singular sextic
hypersurfaces (EPW-sextics) in $\PP^5$ (Example~(9.3) of~\cite{epw})
which come provided with a natural double cover: we have
shown~\cite{og2} that the generic such double cover is a deformation of
the Hilbert square of a $K3$ and that the family of double EPW-sextics
is a locally versal family of projective deformations of $(K3)^{[2]}$.
Thus the family of double EPW-sextics is similar to the family of Fano
varieties of lines on a cubic $4$-fold  (see~\cite{beaudon}), with the
following difference: the Pl\"ucker ample divisor  on the Fano variety
of lines has square $6$ for the Beauville-Bogomolov quadratic form
(see~\cite{beau,beaudon}) while the natural polarization of a double
EPW-sextic has square $2$ (see~\cite{og2}). Let $Y\subset\PP^5$ be a
generic EPW-sextic: we proved in~\cite{og2} that the dual
$Y^{\vee}\subset(\PP^5)^{\vee}$   is another generic EPW-sextic. Thus
we may associate to the natural double cover $X$ of $Y$ a \lq\lq
dual\rq\rq variety $X^{\vee}$ namely the natural double cover of
$Y^{\vee}$. This construction defines a (rational) involution on the
moduli space of double EPW-sextics. In~\cite{og2} we showed that a
generic EPW-sextic is not self-dual and hence the involution on the
moduli space of double EPW-sextics is not the identity; in this paper
we determine the relation between the periods of a double EPW-sextic
and its dual. Before stating the result we recall the definition of
EPW-sextics. Let $V$ be a $6$-dimensional $\CC$-vector space. Choose an
isomorphism $vol\colon\wedge^6 V\overset{\sim}{\lra}\CC$ and let
$\omega$ be the symplectic form on $\wedge^3 V$ defined by wedge
product followed by $vol$. Let $\PP(V)$ be the projective space of
$1$-dimensional sub vector spaces of $V$; then $\omega$ gives $\wedge^3
V\otimes\cO_{\PP(V)}$ the structure of a symplectic vector-bundle of
rank $20$. Let $F$ be  the sub-vector-bundle of $\wedge^3
V\otimes\cO_{\PP(V)}$ whose fiber $F_{[v]}$ over $[v]\in\PP(V)$
consists of tensors divisible by $v$:
\begin{equation}
F_{[v]}:=\{v\wedge w|\ w\in\wedge^2 V\}.
\end{equation}
As is easily checked $F$ is a Lagrangian sub-bundle of $\wedge^3
V\otimes\cO_{\PP(V)}$. Let $\mathbb{LG}(\wedge^3 V)$ be the symplectic
Grassmannian parametrizing $\omega$-Lagrangian subspaces of $\wedge^3
V$. For $A\in \mathbb{LG}(\wedge^3 V)$ we let $\lambda_A$ be the
composition
\begin{equation}\label{mappaquot}
   F \lra \wedge^3
  V\otimes\cO_{\PP(V)}\lra (\wedge^3 V/A)\otimes\cO_{\PP(V)}
\end{equation}
and $Y_A\subset\PP(V)$ be the zero-scheme of $\det(\lambda_A)$; unless
$\det(\lambda_A)$ is identically zero\footnote{If $A$ is the fiber of $F$ over a fixed point then $\det(\lambda_A)$ is identically zero.} $Y_A$
is a sextic hypersurface because
 $\det F\cong\cO_{\PP(V)}(-6)$. An {\it EPW-sextic} is a hypersurface
in $\PP(V)$ which is equal to $Y_A$ for some $A\in \mathbb{LG}(\wedge^3
V)$. In~\cite{og2} we described explicitly the non-empty Zariski-open
$\mathbb{LG}(\wedge^3 V)^0\subset\mathbb{LG}(\wedge^3 V)$ parametrizing
$A$ such that $Y_A$ has no singularities other than those forced by its
description as a degeneracy locus; for $A\in\mathbb{LG}(\wedge^3 V)^0$
the singular locus of $Y_A$ is a smooth degree-$40$ irreducible surface
and at a singular point $Y_A$ is locally (in the classical topology)
isomorphic to the product of $sing(Y_A)$ and an $A_1$-singularity. If
$A\in\mathbb{LG}(\wedge^3 V)^0$ there is a natural double cover
$f_A\colon X_A\to Y_A$ ramified only over $sing(Y_A)$, with $X_A$
smooth.  As shown in~\cite{og2} the $4$-fold $X_A$ is a deformation of
$(K3)^{[2]}$ - the Hilbert square of a $K3$. Let $L_A$ be the
line-bundle on $X_A$ defined by
\begin{equation}\label{ellea}
L_A:=f_A^{*}\cO_{Y_A}(1);
\end{equation}
then $c_1(L_A)$ has square $2$ for the Beauville-Bogomolov quadratic
form. The family of $(X_A,L_A)$ for $A$ varying in
$\mathbb{LG}(\wedge^3 V)^0$ is a locally complete family of polarized
deformations of $(K3)^{[2]}$. Now we recall the duality map. Let
$vol^{\vee}$ be a trivialization of $\wedge^6 V^{\vee}$ and
$\omega^{\vee}$ be the symplectic form on $\wedge^3 V^{\vee}$ given by
wedge-product followed by $vol^{\vee}$: let $\mathbb{LG}(\wedge^3
V^{\vee})$ be the symplectic Grassmannian parametrizing
$\omega^{\vee}$-Lagrangian subspaces of $\wedge^3 V^{\vee}$. For
$A\in\mathbb{LG}(\wedge^3 V)$ the annihilator $A^{\bot}\subset \wedge^3
V^{\vee}$ belongs to $\mathbb{LG}(\wedge^3 V^{\vee})$. Thus we have an
isomorphism
\begin{equation}\label{mappadualita}
  \begin{matrix}
 \delta\colon\mathbb{LG}(\wedge^3 V) &
 \overset{\sim}{\lra} &
 \mathbb{LG}(\wedge^3 V^{\vee})\\
 A & \mapsto & A^{\bot}.
  \end{matrix}
\end{equation}
Let
\begin{equation}
\mathbb{LG}(\wedge^3
 V)^{00}:=\mathbb{LG}(\wedge^3
 V)^0\cap\delta^{-1}\mathbb{LG}(\wedge^3
 V^{\vee})^0.
\end{equation}
Thus $\mathbb{LG}(\wedge^3
 V)^{00}$ is open dense in $\mathbb{LG}(\wedge^3
 V)$.  Assume that $A\in\mathbb{LG}(\wedge^3
 V)^{00}$.
 By Proposition~(3.1)
of~\cite{og2} $Y_{A^{\bot}}$ is $Y^{\vee}_A$ (the dual  of $Y_A$); we
let
\begin{equation}
    X_A^{\vee}\to Y^{\vee}_{A}=Y_{A^{\bot}}
\end{equation}
be the natural double cover and we call $X_A^{\vee}$ the {\it dual of
$X_A$}. We will show how to obtain the Hodge structure on
$H^2(X_A^{\vee})$ from the Hodge structure on $H^2(X_A)$.
 First we describe the relevant period space.  Let
$\wt{\Lambda}$ be the even lattice
\begin{equation}\label{lambdatilde}
 \wt{\Lambda}:=U^3\widehat{\oplus}(-E_8)^2\widehat{\oplus}(-2)
\end{equation}
where $U$ is the hyperbolic plane and  $(-2)$ is the lattice generated
by a single element of square $-2$; we will denote by $(\cdot,\cdot)$
the symmetric bilinear form on $\wt{\Lambda}$. Let $U$ be one of the
hyperbolic planes appearing in Decomposition~(\ref{lambdatilde}) and
choose
\begin{equation}\label{scelgou}
\text{$u\in U$  of square $2$.}
\end{equation}
Let $e_1\in U$ be a generator of $u^{\bot}\cap
U$ and $e_2$ be a generator of the direct summand $(-2)$ appearing in
Decomposition~(\ref{lambdatilde}). Then
 \begin{equation}\label{lambdazero}
\Lambda:=u^{\bot}= U^2\widehat{\oplus}(-E_8)^2\widehat{\oplus}\ZZ
e_1\widehat\oplus\ZZ e_2
 \cong U^2\widehat{\oplus}(-E_8)^2\widehat{\oplus}(-2)^2.
\end{equation}
The period domain of interest to us is
 \begin{equation}
\cD_2:=\{[v]\in\PP(\Lambda\otimes\CC)|\ (v,v)=0,\ \
(v,\overline{v})>0\}.
\end{equation}
Let  $Stab(u)<O(\wt{\Lambda})$ be the subgroup of isometries fixing $u$
and $\rho\colon Stab(u) \to  O(\Lambda)$ be the restriction map - this
makes sense because $\Lambda=u^{\bot}$. Let
\begin{equation}\label{eccogamma}
\Gamma:=Im(\rho)<O(\Lambda).
\end{equation}
The {\it period moduli space} is the quotient
\begin{equation}
 \cQ_2:=\Gamma\backslash \cD_2.
\end{equation}
 We will prove in
Subsection~(\ref{demoindice}) that
\begin{equation}\label{indicedue}
[O(\Lambda):\Gamma]=2.
\end{equation}
An explicit  $r\in(O(\Lambda)\setminus\Gamma)$
is the isometry which interchanges $e_1$ with $e_2$ and fixes the
elements of $\{e_1,e_2\}^{\bot}$, see~(\ref{alfaexpl}).
By~(\ref{indicedue}) the involution $r$ descends to a non-trivial
involution
\begin{equation}\label{etabeta}
\ov{r}\colon\cQ_2\to\cQ_2.
\end{equation}
Let $A\in\mathbb{LG}(\wedge^3 V)^0$: the isomorphism class of the Hodge
structure on $c_1(L_A)^{\bot}\subset H^2(X_A)$ is a point of $\cQ_2$
(see Section~(\ref{mappaper})) that we denote by $\cP(A)$ - this is the
{\it period point of $X_A$}. The global period map is defined by
\begin{equation}\label{perglob}
\begin{matrix}
\mathbb{LG}(\wedge^3 V)^0  & \overset{\cP}{\lra} & \cQ_2 \\
A & \mapsto & \cP(A).
\end{matrix}
\end{equation}
Of course we also have a global period map $\mathbb{LG}(\wedge^3
V^{\vee})^0  \to \cQ_2$; we denote it by the same symbol $\cP$. If
$A\in\mathbb{LG}(\wedge^3 V)^{00}$ we have two period points namely
$\cP(A)$ and $\cP(A^{\bot})$. The following is the main result of this
paper.
\begin{thm}\label{princrisul}
 Keep notation as above.
 If $A\in\mathbb{LG}(\wedge^3 V)^{00}$ then
 \begin{equation}\label{reltraper}
\cP(A^{\bot})=\ov{r}\circ\cP(A).
\end{equation}
 \end{thm}
The theorem may be stated informally as follows: the polarized Hodge structures  on $c_1(L_{A^{\bot}})^{\bot}$ and $c_1(L_{A})^{\bot}$ are always isomorphic while in general there is no Hodge isometry between $H^2(X_A)$ and $H^2(X_{A^{\bot}})$.

There are three steps in the proof of the above theorem. First
of all there exists a \underline{rational} Hodge isometry
\begin{equation}\label{rathodgeiso}
H^2(X_A)\supset c_1(L_A)^{\bot}\cong
c_1(L_{A^{\bot}})^{\bot}\subset H^2(X_{A^{\bot}}).
\end{equation}
In fact let $\wt{Y}_A,\wt{Y}_{A^{\bot}}$ be the desingularizations of $Y_A,Y_{A^{\bot}}$ respectively - they are obtained by blowing up $sing(Y_A)$ and $sing(Y_{A^{\bot}})$ respectively. From the equality $Y_A^{\vee}=Y_{A^{\bot}}$ one gets that
 \begin{equation}\label{pescetti}
 \wt{Y}_A\cong\wt{Y}_{A^{\bot}}.
\end{equation}
Now $H^4(Y_A)$ is a subgroup of $Sym^2 H^2(X_A)$ and similarly for
$H^4(Y_{A^{\bot}})$; this fact together with~(\ref{pescetti})
gives~(\ref{rathodgeiso})  after some work. It follows
from~(\ref{rathodgeiso}) that locally (in the classical topology) there
exists $g\in O(\Lambda\otimes\QQ)$ relating the periods of $X_A$ and
$X_{A^{\bot}}$; this means that we have
\begin{equation}\label{relocale}
\cP_{\Phi}(A^{\bot})=g\circ\cP_{\Psi}(A)
\end{equation}
 where $\cP_{\Psi},\cP_{\Phi}$ are local liftings of the global period
 maps ($\Psi,\Phi$ are markings of the weight-$2$ cohomology).
 The second step consists in showing that~(\ref{reltraper}) holds for $A$ belonging to a certain locally closed  codimension-$1$ submanifold 
 $\Delta^0_{*}(V)\subset\mathbb{LG}(\wedge^3 V)^0$.
If $A\in\Delta^0_{*}(V)$ then $X_A$ is isomorphic to
a moduli space of sheaves on a $K3$ surface $S\subset\PP^6$ of degree
$10$. On the other hand there is a map $f\colon S^{[2]}\to  Y_{A^{\bot}}$ of degree $2$ which contracts a certain plane $P\subset S^{[2]}$ and is finite over $Y_{A^{\bot}}\setminus f(P)$. One shows that $mult_{f(P)}Y_{A^{\bot}}=3$ and hence $A^{\bot}\notin \mathbb{LG}(\wedge^3 V^{\vee})^0$ i.e.~$A\notin \mathbb{LG}(\wedge^3 V)^{00}$. 
Morally $X_{A^{\bot}}$ 
 is the singular symplectic $4$-fold obtained from $S^{[2]}$ 
by contracting the plane $P$. We do not prove a precise version of this; however we do prove that the local period map of $\mathbb{LG}(\wedge^3 V^{\vee})^0$ extends across $A^{\bot}$ and 
that its value  is given by the periods of 
$S^{[2]}$. This will allow us to check that~(\ref{reltraper})
holds for $A\in\Delta^0_{*}(V)$.
Since the local period map sends
$\Delta^0_{*}(V)$ to the intersection of $\cD_2$ and a hyperplane
$\zeta^{\bot}$ it will follow that the rational isometry $g$
of~(\ref{relocale}) is (with suitable markings $\Psi,\Phi$)  either $r$
or $r\circ r_{\zeta}$ where $r_{\zeta}$ is the reflection in the span of $\zeta$. We rule out  this second possibility by a
monodromy argument: this is the third step in the proof of
Theorem~(\ref{princrisul}).

The paper is organized as follows. In Section~(\ref{mappaper}) we
recall the definition of  local and global period maps and we
prove~(\ref{indicedue}). In the next section we prove  the results  on
$\Delta^0_{*}(V)$  that we described
above.  In the final section we prove the existence of a rational Hodge
isometry~(\ref{rathodgeiso}) and we give the monodromy argument.
\vskip 2mm
  \noindent
  {\bf Notation:} $V$ will always be a complex vector-space of
  dimension $6$.
 \section{The period map and the lattice $\Lambda$}\label{mappaper}
\setcounter{equation}{0}
 \subsection{The period map}\label{mappaperiodi}
\setcounter{equation}{0}
Let $X$ be a deformation of $(K3)^{[2]}$. Let  $(\cdot,\cdot)_{X}$ be
the Beauville-Bogomolov symmetric bilinear form on $H^2(X)$.
The restriction of $(\cdot,\cdot)_{X}$ to  $H^2(X;\ZZ)$ is  a
non-degenerate integral symmetric form; thus $H^2(X;\ZZ)$ is a lattice.
As is well-known $H^2(X;\ZZ)$ and $\wt{\Lambda}$ are
isometric. Now assume that we are
given a holomorphic line-bundle $L$ on $X$ such that
$(c_1(L),c_1(L))_{X}=2$. Since $O(\wt{\Lambda})$ acts transitively on the set of
vectors of square $2$  there exists an isometry
\begin{equation}\label{marcato}
 \psi\colon H^2(X;\ZZ)\overset{\sim}\lra\wt{\Lambda}
\end{equation}
such that $\psi(c_1(L))=u$ where $u$ is as in~(\ref{scelgou});
this is a {\it marking of $(X,L)$}.  Let $\psi_{\CC}\colon H^2(X;\CC)\overset{\sim}\lra\wt{\Lambda}\otimes\CC$ be the map
obtained from $\psi$ by extension of scalars.
Let $\sigma$ be a generator of $H^{2,0}(X)$; then
\begin{equation}\label{prodaccadue}
(c_1(L),\sigma)_{X}=0=(\sigma,\sigma)_{X},\quad
(\sigma,\ov{\sigma})_{X}>0.
\end{equation}
Thus
\begin{equation}
\cP_{\psi}(X,L):=\psi_{\CC}(H^{2,0}(X))\in\cD_2.
\end{equation}
This is the {\it local period point of $(X,L)$ associated to $\psi$}
(or periods of $(X,L)$). Any other marking  of $(X,L)$ is given by
$\gamma\circ\psi$ where $\gamma\in\Gamma$, thus the $\Gamma$-orbit of
$\cP_{\psi}(X,L)$ is a well-defined point of $\cQ_2$; this is {\it the
global period point of $(X,L)$}, we denote it by $\cP(X,L)$. If
$A\in\mathbb{LG}(\wedge^3 V)^0$ we let $\cP(A):=\cP(X_A,L_A)$.

Later on we will study global and local  period maps for certain
families of double EPW-sextics. Let
\begin{equation}\label{famiglia}
\pi\colon\cX\to T
\end{equation}
be a proper submersive map between complex manifolds such that each
fiber is a deformation of $(K3)^{[2]}$. For $t\in T$ we let
$X_t:=\pi^{-1}(t)$. We assume that we are given a (holomorphic)
line-bundle $\cL$ on $\cX$ such that $c_1(\cL|_{X_t})$ has square $2$
for every $t\in T$; we let $L_t:=\cL|_{X_t}$. The
{\it global period map of $(\cX,\cL)$} is given by
\begin{equation}
\begin{matrix}
T & \overset{\cP}{\lra} & \cQ_2 \\
t & \mapsto & \cP(t):=\cP(X_t,L_t)
\end{matrix}
\end{equation}
Griffiths proved that  $\cP$ is a holomorphic map. Let $\pi\colon\cX\to
\mathbb{LG}(\wedge^3 V)^0$ be the tautological family of double
EPW-sextics  and $\cL$ be the tautological relatively ample line-bundle
on $\cX$ which restricts to $L_A$ (see~(\ref{ellea})) on $X_A$; then
$\cP$
 is  the period map of~(\ref{perglob}), in particular the map
 of~(\ref{perglob})
 is holomorphic. Let's go back to a general family~(\ref{famiglia}): it
 is not always possible to lift the global period map $\cP$ to a map
 $T\to\cD_2$, in fact a necessary and sufficient condition is that
$R^2\pi_{*}\ZZ$ is trivial. Suppose that $R^2\pi_{*}\ZZ$ is trivial.
Then there exists a trivialization $\Psi\colon
R^2\pi_{*}\ZZ\overset{\sim}{\lra}T\times\wt{\Lambda}$ sending the
section corresponding to $c_1(\cL)$ to the section given by $u$ - this
is a {\it marking of $(\cX,\cL)$}. Given such a marking  we let
$\Psi_t\colon H^2(X_t;\ZZ)\overset{\sim}{\lra}\wt{\Lambda}$ be the
fiber of $\Psi$ over $t$. The {\it local period map of $(\cX,\cL)$
associated to $\Psi$} is given by
\begin{equation}
\begin{matrix}
T & \overset{\cP_{\Psi}}{\lra} & \cD_2 \\
t & \mapsto & \cP_{\Psi}(t):=\cP_{\Psi_t}(X_t,L_t).
\end{matrix}
\end{equation}
Griffiths proved that $\cP_{\Psi}$ is holomorphic.
 \subsection{Proof of~(\ref{indicedue})}\label{demoindice}
\setcounter{equation}{0}
We start by recalling  the definition  of discriminant group and
discriminant quadratic form of an even lattice $(L,(\cdot,\cdot)_L)$,
i.e.~a free finitely generated abelian group
 $L$ equipped with a symmetric integral  even non-degenerate bilinear form
 $(\cdot,\cdot)_L$. We follow~\cite{nik}. The bilinear form
 $(\cdot,\cdot)_L$ extends
 to a $\QQ$-valued bilinear form on $L\otimes\QQ$; abusing notation we denote
by $(\cdot,\cdot)_L$ the extended form. Let $L^{\vee}:=Hom(L,\ZZ)$; by
non-degeneracy of $(\cdot,\cdot)_L$ we have  a natural chain of
inclusions
\begin{equation}\label{tuttidentro}
    L\subset L^{\vee}\subset L\otimes\QQ.
\end{equation}
The {\it discriminant group of $L$} is $A_L:=L^{\vee}/L$; it comes
provided with the {\it discriminant bilinear-form}
\begin{equation}
\begin{matrix}
A_L\times A_L & \overset{b_L}{\lra} & \QQ/\ZZ  \\
([\alpha],[\beta]) & \mapsto & [(\alpha,\beta)_L]
\end{matrix}
\end{equation}
and the {\it discriminant quadratic-form}
\begin{equation}
\begin{matrix}
A_L & \overset{q_L}{\lra} & \QQ/2\ZZ  \\
[\alpha] & \mapsto & [(\alpha,\alpha)_L].
\end{matrix}
\end{equation}
The formula
\begin{equation}
q_L([\alpha+\beta])\equiv q_L([\alpha])+q_L([\beta])+2
b_L([\alpha],[\beta]) \pmod{2\ZZ}
\end{equation}
shows that $q_L$ determines uniquely $b_L$. An {\it index-$i$
overlattice of $L$} consists of a lattice $M$ and an inclusion of
lattices $L\subset M$ (the restriction of $(\cdot,\cdot)_M$ to $L$ is
equal to $(\cdot,\cdot)_L$) of  index $i$. Two overlattices $M_1\supset
L\subset M_2$ are equivalent if there exists an isometry
$M_1\overset{\sim}{\lra} M_2$ which is the identity on $L$, i.e.~if the
inclusions $M_i\hra L^{\vee}$ for $i=1,2$ have the same image. Suppose
that $L\subset M$ is an  index-$i$ overlattice of $L$; the inclusion
$M\subset L^{\vee}$ defines an inclusion $M/L\subset A_L$ with image a
subgroup of cardinality $i$ which is $q_L$-isotropic. This construction
defines a one-to-one correspondence between the set of equivalence
classes of indexi-$i$ overlattices of $L$  and the set of
$q_L$-isotropic subgroups of $A_L$ of cardinality $i$. The
correspondence is equivariant for the natural actions of $O(L)$ on both
sets. Now consider the lattice $L=\ZZ u\widehat{\oplus}\Lambda$ where $u$ and $\Lambda$ are as in
 Section~(\ref{prologo}):
 we will describe the discriminant group
and discriminant form of $L$.  Let $e_1,e_2\in\Lambda$ be as in
Section~(\ref{prologo}). A straightforward computation gives that
\begin{equation}
    A_L=\ZZ/(2)\left[\frac{u}{2}\right]\oplus
    \ZZ/(2)\left[\frac{e_1}{2}\right]
    \oplus\ZZ/(2)\left[\frac{e_2}{2}\right]
\end{equation}
 and that
 \begin{equation}
 q_L\left(x\left[\frac{u}{2}\right]+y_1\left[\frac{e_1}{2}\right]+
 y_2\left[\frac{e_2}{2}\right]\right)\equiv
 \frac{1}{2}(x^2-y_1^2-y_2^2)\pmod{2\ZZ}.
\end{equation}
The set $I\subset A_L$ of non-zero isotropic vectors is given by
\begin{equation}
I=\left\{\left[\frac{u}{2}\right]+\left[\frac{e_1}{2}\right],
\left[\frac{u}{2}\right]+\left[\frac{e_2}{2}\right]\right\}.
\end{equation}
The group $O(L)$ acts naturally on $A_L$ and hence also on $I$; thus we
have a homomorphism
\begin{equation}
O(L)\overset{\epsilon}{\lra} Aut(I)\cong\ZZ/(2).
\end{equation}
The overlattice $\wt{\Lambda}\supset L$ is of index $2$ because $L$ is
the kernel of the surjection
\begin{equation}
\begin{matrix}
 \wt{\Lambda} & \lra & \ZZ/(2) \\
 v & \mapsto & (v,u)\pmod{2}.
\end{matrix}
\end{equation}
The correspondence described above  defines an $O(L)$-equivariant
bijective map between $I$ and the set of equivalence classes of
index-$2$ overlattices of $L$; thus
\begin{equation}\label{nucleoeps}
    Im(O(\wt{\Lambda})\to O(L))=\ker(\epsilon).
\end{equation}
The subgroup of $O(L)$ consisting of isometries which are the identity
on $\ZZ u$ is naturally identified with $O(\Lambda)$; thus $O(\Lambda)<
O(L)$. By~(\ref{nucleoeps}) and the definition of $\Gamma$
(see~(\ref{eccogamma})) we get that
\begin{equation}\label{gammanuc}
\Gamma=\ker(\epsilon|_{O(\Lambda)}).
\end{equation}
Let $r\in O(\Lambda)$ be the involution characterized by the following
properties:
\begin{equation}\label{alfaexpl}
 r(e_1)=e_2, \quad r(e_2)=e_1,\quad
 r|_{\{e_1,e_2\}^{\bot}}=\text{identity}.
\end{equation}
Then $\epsilon(r)$ is the non-trivial permutation of $I$ and hence
$\epsilon|_{O(\Lambda)}$ is surjective. Thus~(\ref{indicedue}) follows
from~(\ref{gammanuc}).
 \section{Explicit dual couples}\label{coppieduali}
\setcounter{equation}{0}
Let $F^3_5\subset\PP^6$ be the intersection of
$Gr(2,\CC^5)\subset\PP(\wedge^2 \CC^5)\cong\PP^9$ with a transversal
$6$-dimensional linear subspace of $\PP^9$. Let
$T\subset|\cO_{\PP^6}(2)|$ be the open dense subset parametrizing
quadrics which are transversal to $F^3_5$. For $t\in T$ let $Q_t$ be
the quadric corresponding to $t$ and
\begin{equation}\label{kappatre}
  S_t:=F^3_5\cap Q_t.
\end{equation}
Then $S_t$ is a degree-10 linearly normal $K3$ surface; in fact the
generic such surface is projectively equivalent to $S_t$ for some $t\in
T$ by Mukai~\cite{mukaik3}. Let $D_t$ be the hyperplane divisor class
on $S_t$ and $M_t$ be the moduli space of $D_t$-semistable rank-$2$
sheaves $\cF$ on $S_t$ with $c_1(\cF)=c_1(D_t)$ and $c_2(\cF)=5$.
Suppose that
\begin{equation}\label{dieci}
  \text{ for all divisors $E$  on $S_t$ we have $E\cdot D_t\equiv 0\pmod{10}$.}
\end{equation}
Then (see Section~(5) of~\cite{og2}) there exists $A_t \in
\mathbb{LG}(\wedge^3 V)^0$ such that
\begin{equation}\label{emmeti}
X_{A_t}\cong M_t.
\end{equation}
Furthermore $Y_{A_t}$ is explicitely described as follows. Let
$\Sigma_t$ be the divisor on $|I_{S_t}(2)|$ parametrizing singular
quadrics, since all quadrics parametrized by $|I_{F^3_5}(2)|$ are
singular we have
\begin{equation}\label{discrimine}
\Sigma_t=|I_{F^3_5}(2)|+\Sigma'_t.
\end{equation}
Then
\begin{equation}\label{ipsilonti}
Y_{A_t}\cong \Sigma'_t.
\end{equation}
(Formally $\Sigma'_t$ is a degree-$6$ divisor; in the above equation we are implicitely stating that $\Sigma'_t$ is a reduced divisor and hence we may view it as a degree-$6$ hypersurface.)
The set of $t$ for which~(\ref{dieci}) holds is the  complement of a
countable union of proper algebraic subvarieties of $T$ however a
straightforward argument shows that there is an open dense $T'''\subset
T$ such that for $t\in T'''$ there exists $A_t \in \mathbb{LG}(\wedge^3
V)^0$ for which both~(\ref{emmeti}) and~(\ref{ipsilonti}) hold; we give
the argument in Subsection~(\ref{deltazero}).  We let
$\Delta^0(V)\subset\mathbb{LG}(\wedge^3 V)^0$ be the set of $A$ such
that $X_A$ is isomorphic to $M_t$ for some $t\in T'''$; this is a
locally closed subset of $\mathbb{LG}(\wedge^3 V)^0$. Computing the
periods of $X_A$ for $A\in\Delta^0(V)$ we will show that $\Delta^0(V)$
has codimension $1$ in $\mathbb{LG}(\wedge^3 V)^0$. Let
$A\in\Delta^0(V)$, thus $X_A\cong M_t$ for some $t\in T$: the dual
$Y_{A}^{\vee}$ is described as follows. The Hilbert scheme $S_t^{[2]}$
contains a copy of $\PP^2$, call it $P_t$, parametrizing $Z\subset S_t$ which span a line contained in $F^3_5$. Let $S_t^{[2]}\to N_t$ be the contraction of  $P_t$ - thus $N_t$ is a singular symplectic variety.
There is an involution on $N_t$ whose quotient is isomorphic to
$Y_{A}^{\vee}$. From this it will follow that
$A\notin\mathbb{LG}(\wedge^3 V)^{00}$ and hence  apparently it will not
make sense to compute $\cP(A^{\bot})$. The main point of this
subsection is to prove that the  period map extends across $A^{\bot}$,
in fact the \underline{local} period map extends and its value at
$A^{\bot}$ is given by the periods of $(S_t^{[2]},D_t)$.
\subsection{The locus $\Delta^0(V)$}
\label{deltazero}
\setcounter{equation}{0}
We recall that the {\it Mukai vector $v(\cF)$} of a sheaf $\cF$ on
$S_t$ is
\begin{equation}\label{vettoremukai}
v(\cF):=ch(\cF)\sqrt{Td(S_t)}=ch(\cF)(1+\eta_t)\in H^{*}(S_t;\ZZ),
\end{equation}
where $\eta_t\in H^4(S_t;\ZZ)$ is the orientation class. If $[\cF]\in
M_t$ (we let $[\cF]\in M_t$ be the  equivalence class of the semistable
sheaf $\cF$ - we recall that if $\cF$ is stable this is the same as the
isomorphism class of $\cF$) the class $v(\cF)$ is independent of $\cF$,
we denote it by $v_t$; explicitly
\begin{equation}\label{vuti}
v_t=2+c_1(D_t)+2\eta_t.
\end{equation}
\begin{prp}\label{emmetiepw}
Keep notation as above. There is an open dense subset $T'''\subset T$
such that the following holds. Let $t\in T'''$; then there exists
$A_t\in\mathbb{LG}(\wedge^3 V)^{0}$ such that both~(\ref{emmeti})
and~(\ref{ipsilonti}) hold. In particular we have a canonical
identification $\PP(V)\cong | I_{S_t}(2)|$.
\end{prp}
\begin{proof}
By Maruyama~\cite{maru} there exists a projective map $\rho\colon\cM\to
T$ with (schematic) fiber $M_t$ over $t\in T$. Let $t\in T$; we say
that $S_t$ is {\it unsuitable} if there exists a divisor class $C$ on
$S_t$ such that
\begin{equation}
\int_{S_t} c_1(C)\wedge c_1(D_t)=0,\quad
 -10\le \int_{S_t} c_1(C)^2 < 0.
\end{equation}
The set of unsuitable $t$ is a proper closed subset of $T$, thus the
complement  $T'$  is an open dense subset of $T$. Let
$\cM':=\rho^{-1}(T')$ and $\rho'\colon\cM'\to T'$ be the restriction of
$\rho$. It is known that if $t\in T'$ then every sheaf parametrized by
$M_t$ is slope-stable and $M_t$ is a  smooth $4$-dimensional scheme
(Main Theorem~(0.1.2) and Proposition~(2.1) of~\cite{ogradyvb} - notice
that $t\in T_1$ if and only if $D_t$ is $|v_t|$-generic). Let's show
that there is an open dense  $T''\subset T'$ such that for $t\in T''$
the Mukai reflection (see (4.2.2) of~\cite{og1}) is a regular
involution $\phi_t$ on $M_t$ - we recall that in general the Mukai
reflection acts on the derived category of coherent sheaves on $S_t$.
In order for $\phi_t$ to be a regular involution on $M_t$ it suffices
that for all $[\cF]\in M_t$
 the following hold:
\begin{itemize}
\item [(a)]
$h^0(\cF)=\chi(\cF)=4$,
\item [(b)]
$\cF$ is globally generated away from (at most) a zero-dimensional
subset of $S_t$,
\item [(c)]
the kernel of the evaluation map $H^0(\cF)\otimes\cO_{S_t}\to \cF$,
call it $\cE$, is a $D_t$-slope-stable sheaf.
\end{itemize}
If~(a)-(b)-(c) are satisfied for all $[\cF]\in M_t$ then the generic
sheaf $\cF$ parametrized by $M_t$  is globally generated and for such a
sheaf $\phi_t([\cF])=[\cE^{\vee}]$ where $\cE$ is as in~(c) above.
Let $T''\subset T'$ be the set of $t$ such that~(a)-(b)-(c) hold for every $[\cF]\in M_t$. Let's prove that $T''$ is open.
First we show that the subset
 $T_a\subset T$  of $t\in T'$ such that~(a) holds for every $[\cF]\in M_t$ is open.  Let $t\in T'$: if
$[\cF]\in M_t$ then by stability $h^2(\cF)=0$ and hence $h^0(\cF)\ge
4$, thus by the upper-semicontinuity Theorem the set of $[\cF]\in M_t$ that violate~(a) is  closed.  By properness of the map $\rho'\colon\cM'\to T'$ it follows that $(T'\setminus T_a)$  is closed i.e.~$T_a$ is open. A similar argument shows that the set of $t\in T_a$ such that~(b)-(c) hold for every $[\cF]\in M_t$ is open;
thus $T''$ is open.
Since $T''$ contains the subset of $t$ such
that~(\ref{dieci}) holds it is dense in $T$. Let $T'''\subset T''$ be
the set of $t$ such that
\begin{equation}\label{retteconiche}
\text{$S_t$ contains no effective non-zero divisor $E$ with
$E\cdot D_t\le 5$.}
\end{equation}
Thus $T'''$ is open and dense in $T$. We claim that if $t\in T'''$ then
\begin{equation}\label{emmesigma}
M_t/\langle\phi_t
\rangle\cong \Sigma'_t\subset |I_{S_t}(2)|\cong\PP^5.
\end{equation}
More precisely Proposition~(5.1) of~\cite{og1} holds for
$M_t=\cM(v_t)$. In order to prove this it suffices to show that for
every $t\in T'''$ the following holds:
\begin{itemize}
\item[(1)]
If $[\cF]\in M_t$ and $\sigma\in H^0(\cF)={\rm Hom}(\cO_{S_t},\cF)$ is
non-zero then the quotient $\cF/Im(\sigma)$ is locally-free in
codimension $1$ ($\sigma$ has isolated zeroes) - this is used in the
proof of Lemma~(5.4) of~\cite{og1}.
\item[(2)]
If $\cG$ is a globally generated rank-$2$ vector bundle on $S_t$ with
$\det\cG\cong\cO_{S_t}(D_t)$ and $c_2(\cG)=5$ then $\cG$ is $D_t$-slope
stable - this is used in the proof of Lemma~(5.7) of~\cite{og1}.
\end{itemize}
Let's show that~(1) above holds. Suppose that $\sigma$ does not have
isolated zeroes: then it vanishes along a non-zero effective divisor
$E$ and by slope stability of $\cF$ we have $E\cdot D_t< slope(\cF)=5$
contradicting~(\ref{retteconiche}).  Let's show that~(2) above holds.
Suppose that $\cG$ is not $D_t$-slope stable; since $\cG$ has rank $2$
there is a destabilizing  sequence $\cG\to I_Z\otimes\cO_{S_t}(E)$
where $Z\subset S_t$ is zero-dimensional and $E$ is a divisor with
$E\cdot D_t\le slope(\cF)=5$. Since $\cG$ is globally generated there
is a non-zero section of  $I_Z\otimes\cO_{S_t}(E)$ and hence $E$ is
effective; by~(\ref{retteconiche}) we get that $E=0$. Thus we have an
exact sequence
\begin{equation}
0\lra \cO_{S_t}(D_t)\lra \cG\lra I_Z \lra 0.
\end{equation}
Since $\ell(Z)=c_2(\cG)=5$ the zero-dimensional subscheme $Z\subset
S_t$ is non-empty and hence $h^0(I_Z)=0$; this contradicts the
hypothesis that $\cG$ is globally generated. We have proved
that~(\ref{emmesigma}) holds for $t\in T'''$. By Theorem~(1.1)
of~\cite{og2} we get that there exists $A_t\in\mathbb{LG}(\wedge^3
V)^0$ such that both~(\ref{emmeti}) and~(\ref{ipsilonti}) hold.
\end{proof}
The proof of the above proposition together with Claim~(5.18)
of~\cite{og1} gives the following result.
\begin{crl}\label{effephi}
Let $t\in T'''$ and let $A_t\in\mathbb{LG}(\wedge^3 V)^0$ be such that
both~(\ref{emmeti}) and~(\ref{ipsilonti}) hold. The map $f_{A_t}\colon
X_{A_t}\to Y_{A_t}$ is identified with the quotient map $M_t\to
M_t/\langle\phi_t\rangle$.
\end{crl}
Let $\Delta^0(V)\subset \mathbb{LG}(\wedge^3 V)^0$ be the locus of $A$
such that $X_A$ is isomorphic to $M_t$ for some $t\in T'''$. Our next
task is to show that $\Delta^0(V)$ is a locally closed subset of
codimension $1$. First we recall how one describes $H^2(M_t)$  for
$t\in T'$. Let $u,w\in H^{*}(S_t)$ and let $u_q,w_q$ be the degree-$q$
components of $u,w$ respectively. One sets $u^{\vee}:=u_0-u_2+u_4$
 and
\begin{equation}\label{mukform}
\langle u,w\rangle:=\int_{S_t} (u_2\wedge w_2-u_0 w_4-u_4 w_0)
=-\int_{S_t}u^{\vee}\wedge w.
\end{equation}
This is {\sl Mukai's bilinear symmetric form\/}.  One defines a
positive\footnote{This means that $h^{p,q}=0$ if $p<0$.} weight-two
Hodge structure on $H^{*}(S_t)$ by defining the Hodge filtration as
\begin{equation}
F^1 H^{*}(S_t):=H^0(S_t)\oplus F^1 H^2(S_t)\oplus H^4(S_t),\quad F^2
H^{*}(S_t):=F^2 H^2(S_t).
\end{equation}
Let $v_t$ be the Mukai vector~(\ref{vuti}),  then $v_t$ is integral of
type $(1,1)$ and hence $v_t^{\bot}$ is an integral Hodge substructure
of $H^{*}(S_t)$ and furthermore the restriction of Mukai's bilinear
symmetric form to $v_t^{\bot}$ is integral. Mukai defined
(see~\cite{mukaisug,ogradyvb}) a map
\begin{equation}\label{thetaiso}
 \theta_t\colon v_t^{\bot}\lra H^2(M_t)
 \end{equation}
by taking K\"unneth components of the Chern character of a tautological
sheaf on $S_t\times M_t$ (if such a sheaf does not exists one considers
a quasi-tautological sheaf).  In~\cite{ogradyvb} we proved that
$\theta_t$ is an isomorphism of Hodge structures and an isometry of
lattices - of course the bilinear form on $v_t^{\bot}$ is the
restriction of the Mukai pairing.  Now assume that $t\in T'''$ and let
$A_t$ be as in Proposition~(\ref{emmetiepw}) and $L_{A_t}$ be as
in~(\ref{ellea}); then by  Corollary~(\ref{effephi}) and
 Proposition~(5.1) of~\cite{og1} we have
 \begin{equation}\label{elleta}
 c_1(L_{A_t})=\theta_t(\eta_t-1).
\end{equation}
We let $h_t:=\theta_t(\eta_t-1)$.
\begin{prp}\label{codelta}
Keep notation as above. Then $\Delta^0(V)$ is a $PGL(V)$-invariant
subset of $\mathbb{LG}(\wedge^3 V)^0$ which is locally (in the
classical topology) a codimension $1$ submanifold.
\end{prp}
\begin{proof}
The subset $\Delta^0(V)$ is  $PGL(V)$-invariant by definition. Let
$\pi\colon\cX\to \mathbb{LG}(\wedge^3 V)^0$ be the tautological family
of double EPW-sextics and $\cL$ be the tatutological relatively ample
line-bundle on $\cX$; thus the restriction of $\cL$ to $X_A$ is
isomorphic to $L_A$. Let $A_p\in \Delta^0(V)$. Thus there exists $p\in
T'''$ such that $X_{A_p}\cong M_p$. Let $U\subset \mathbb{LG}(\wedge^3
V)^0$ be a small open ball contining $A_p$. Let $\cX_U:=\pi^{-1}(U)$
and $\cL_U:=\cL|_{\cX_U}$. Since $U$ is contractible there is a marking
$\Psi$ of $(\cX_U,\cL_U)$; let $\cP_{\Psi}\colon U\to\cD_2$ be the
corresponding  local period map. We notice that
$(5+2c_1(D_p)+5\eta_p)\in v_p^{\bot}$ and hence
$\theta_p(5+2c_1(D_p)+5\eta_p)\in H^2(M_p)$. Furthermore since $\langle
\eta_p-1, 5+2c_1(D_p)+5\eta_p\rangle=0$ we have
$(h_p,\theta_p(5+2c_1(D_p)+5\eta_p))_{M_p}=0$; since $h_p=c_1(L_{A_p})$
we get that $\Psi_p(\theta_p(5+2c_1(D_p)+5\eta_p))\in u^{\bot}$ where
$u$ is our fixed vector of square $2$ - see~(\ref{scelgou}). Let's
prove that
\begin{equation}\label{imagodelta}
\cP_{\Psi}(\Delta^0(V)\cap U)=\cP_{\Psi}(U)\cap
\Psi_p\circ\theta_p(5+2c_1(D_p)+5\eta_p)^{\bot}.
\end{equation}
 First we prove that
\begin{equation}\label{sindentro}
\cP_{\Psi}(\Delta^0(V)\cap U)\subset\cP_{\Psi}(U)\cap
\Psi_p\circ\theta_p(5+2c_1(D_p)+5\eta_p)^{\bot}.
\end{equation}
Let $\cM''':=\rho^{-1}(T''')$ and $\rho'''\colon\cM'''\to T'''$ be the
restriction of $\rho$. Then $\theta_t(5+2c_1(D_t)+5\eta_t)$ is a flat
section of $R^2\rho'''_{*}\ZZ$ and for all $t\in T'''$ we have
\begin{equation}
(\theta_t(5+2c_1(D_t)+5\eta_t),H^{2,0}(S_t))_{M_t}=0
\end{equation}
 because $\theta_t(5+2c_1(D_t)+5\eta_t)\in
H^{1,1}(M_t)$;   the above equality gives~(\ref{sindentro}). Next
we prove that
\begin{equation}\label{destradentro}
\cP_{\Psi}(\Delta^0(V)\cap U)\supset\cP_{\Psi}(U)\cap
\Psi_p\circ\theta_p(5+2\theta_p(c_1(D_p)+5\eta_p)^{\bot}.
\end{equation}
First we notice that
\begin{equation}
\{v_p,\eta_p-1,5+2c_1(D_p)+5\eta_p\}^{\bot}=H^2(S_p)_{prim}
\end{equation}
where the primitive cohomology $H^2(S_p)_{prim}\subset H^2(S)$ is the
orthogonal to $c_1(D_p)$. Since $h_p=\theta_p(\eta_p-1)$ we get that
\begin{equation}\label{accaprim}
H^2(M_p)\supset\{h_p,\theta_p(5+2c_1(D_p)+5\eta_p)\}^{\bot}=
 \theta_p(H^2(S_p)_{prim}).
\end{equation}
Let $\cK_{10}$ be the period space
 for $K3$ surfaces with a polarization of degree $10$;
 Equality~(\ref{accaprim}) defines an isomorphism
\begin{equation}\label{isoper}
 \cD_2\cap\Psi_p\circ\theta_p(5+2c_1(D_p)+5\eta_p)^{\bot}
 \overset{\sim}{\lra}\cK_{10}
\end{equation}
which is compatible with  local period maps defined by the family
$\rho'''\colon \cM'''\to T'$ and the family $\zeta\colon\cS'''\to T'''$
with fiber $S_t$ over $t\in T'''$. Let $\cS'''_U:=\zeta^{-1}(U)$. Since
$\cS'''$ contains the generic $K3$ of degree $10$ the local period map
of the family $\cS'''_U\to U$ is a submersive map from $U$ to an open
ball in $\cK_{10}$; since~(\ref{isoper}) is an isomorphism this
proves~(\ref{destradentro}). We also get that $\cP_{\Psi}$ is
submersive and hence in order to show that $\Delta^0(V)\cap U$ is a
codimension $1$ submanifold it suffices to prove that
\begin{equation}\label{sezione}
\cD_2\cap\Psi_p\circ\theta_p(5+2c_1(D_p)+5\eta_p)^{\bot}
\end{equation}
is smooth. The period domain $\cD_2$ is an open subset of the quadric
of isotropic lines for the non-degenerate quadratic form
$(,)|_{\Lambda\otimes\CC}$ and hence if
$\Psi_p\circ\theta_p(5+2c_1(D_p)+5\eta_p)$ is not isotropic
then~(\ref{sezione}) is smooth. Since
\begin{multline}
\left(\Psi_p\circ\theta_p(5+2c_1(D_p)+5\eta_p),
\Psi_p\circ\theta_p(5+2c_1(D_p)+5\eta_p)\right)= \\
 =\langle 5+2c_1(D_p)+5\eta_p,5+2c_1(D_p)+5\eta_p \rangle=-10
\end{multline}
the intersection~(\ref{sezione}) is indeed smooth.
\end{proof}
\subsection{The dual of $Y_A$ for $A\in\Delta^0(V)$}\label{dualesestica}
\setcounter{equation}{0}
If $A\in \mathbb{LG}(\wedge^3 V)^{00}$ then $Y_{A}^{\vee}=Y_{A^{\bot}}$
by Proposition~(3.1) of~\cite{og2}. As we will see $\Delta^0(V)\cap
\mathbb{LG}(\wedge^3 V)^{00}=\es$ and hence  in order to show that
$Y_{A}^{\vee}=Y_{A^{\bot}}$ for $A\in \Delta^0(V)$ we need to improve
on  the result of~\cite{og2}.
\begin{prp}\label{interduale}
Let $A\in\mathbb{LG}(\wedge^3 V)$ and $\PP(W)\in\PP(V^{\vee})$. Then
$\PP(W)\in Y_{A^{\bot}}$ if and only if
\begin{equation}\label{nonvuoto}
\wedge^3 W\cap A\not=\emptyset.
\end{equation}
\end{prp}
\begin{proof}
Let $\phi\in V^{\vee}$ be a linear function such that $W=\ker(\phi)$:
then
\begin{equation}
F_{\phi}:=\{\phi\wedge\psi|\ \psi\in\wedge^2 V^{\vee}\}=
 \left(\wedge^3 W\right)^{\bot}.
\end{equation}
By definition $\PP(W)\in Y_{A^{\bot}}$ if and only if
\begin{equation}\label{obladi}
\{0\}\not= F_{\phi}\cap A^{\bot}=(\wedge^3 W+A)^{\bot}.
\end{equation}
Since $10=\dim(\wedge^3 W)=\dim A$ and $\dim(\wedge^3 V)=20$ we get
that~(\ref{obladi}) holds if and only if~(\ref{nonvuoto}) holds.
\end{proof}
\begin{crl}\label{nontutto}
Let $A\in\mathbb{LG}(\wedge^3 V)^0$. Then $Y_{A^{\bot}}$ is a
hypersurface and $Y_A^{\vee}=Y_{A^{\bot}}$.
\end{crl}
\begin{proof}
Let $\PP(W)\in Y_{A^{\bot}}$; by Proposition~(\ref{interduale}) this is
equivalent to~(\ref{nonvuoto}). Let $0\not=\alpha\in(\wedge^3 W \cap
A)$: since $\dim W=5$ there exists $v\in W$ such that $\alpha$ is
divisible by $v$ and hence $[v]\in Y_A$. We decompose $W=\CC v\oplus
W_0$ and write $\alpha=v\wedge w$ where $w\in \wedge^2 W_0$; by
Definition~(2.5) of~\cite{og2} the rank of $w$ is $4$, thus
$W=span(v,w)$. Furthermore if $\dim (F_{[v]}\cap A)=1$ then $Y_A$ is
smooth at $[v]$ and $\PP(W)$ is the projective tangent space to $Y_A$
at $[v]$ - see the proof of Proposition~(3.1) of~\cite{og2}. By
Proposition~(2.8) of~\cite{og2} we know that $\dim (F_{[v]}\cap A)=1$
unless $[v]\in sing(Y_A)$ and in this case $\dim (F_{[v]}\cap A)=2$.
Since $sing(Y_A)$ is a surface a straightforward dimension count gives
that $\PP(W)\notin Y_{A^{\bot}}$ for generic $\PP(W)\in\PP(V^{\vee})$,
thus $Y_{A^{\bot}}$ is a hypersurface. The same dimension count gives
that the generic $\PP(W)\in Y_{A^{\bot}}$ is tangent to $Y_A$ at one of
its smooth points; this proves the corollary.
\end{proof}
We will describe explicitely $Y_A^{\vee}=Y_{A^{\bot}}$ for
$A\in\Delta^0(V)$; essentially we will give a refinement of
Proposition~(5.20) and Corollary~(5.21) of~\cite{og1}. Let  $t\in T'''$
and  $S_t$ be  the $K3$  surface corresponding  to $t$; in order to
simplify notation we temporarily drop the subscript $t$. Let $R$ be the
Fano variety of lines on $F^3_5$. If $[\ell]\in R$ then
$\ell\not\subset S$ by~(\ref{retteconiche})
 and hence $\ell\cap Q$ is a
$0$-dimensional scheme of length $2$ contained in $S$: thus we have a
regular map
\begin{equation}\label{mappainter}
\begin{matrix}
R & \lra & S^{[2]} \\
\ell & \mapsto & \ell\cap Q.
\end{matrix}
\end{equation}
Let $P\subset S_t^{[2]}$ be the image of the above map:
then~(\ref{mappainter}) defines a regular map $R\to P$ with inverse
given by
\begin{equation}
\begin{matrix}
P & \lra & R \\
[Z] & \mapsto & span(Z)
\end{matrix}
\end{equation}
and hence $P$ is isomorphic to $R$. It is known~\cite{iskovskih} that
$R\cong\PP^2$, thus $P\cong\PP^2$. Since $S^{[2]}$ is a symplectic
variety it follows that we can contract $P$:
\begin{equation}\label{contraggo}
c\colon S^{[2]}\to N.
\end{equation}
A priori $N$ is a complex space, we will show soon that it is
projective. Let $p:=c(P)$; thus $p$ is the unique singular point of
$N$. On $S^{[2]}$ there is an interesting map to $|I_S(2)|^{\vee}$,
see~(4.3) of~\cite{og1}; we recall the definition. The $K3$ surface $S$
is cut out by quadrics and it contains no lines
by~(\ref{retteconiche}); thus we have a regular map
\begin{equation}\label{mappaeffe}
\begin{matrix}
S^{[2]} & \lra & |I_S(2)|^{\vee}\cong\PP^5\\
[Z] & \mapsto & \{Q_{\lambda}\in |I_S(2)|\ |
\ span(Z)\subset Q_{\lambda}\}.
\end{matrix}
\end{equation}
Let $W\subset |I_S(2)|^{\vee}$ be the image of the above map;
thus~(\ref{mappaeffe}) defines a map $f\colon S^{[2]}\to W$.
 If $[Z]\in P$ then $f([Z])=| I_{F^3_5}(2)|$ hence $f$ is constant on
 $P$; we will see that the point
 \begin{equation}\label{puntospeciale}
f(P)=| I_{F^3_5}(2)|\in W
\end{equation}
is quite special.
 Since $f$ is constant on $P$ and $N$ is normal the
 map $f$ descends to a regular
 map
 \begin{equation}\label{quozdin}
\ov{f}\colon N\to  W.
\end{equation}
\begin{lmm}\label{gradoeffe}
Keep notation as above. There exist a non trivial involution
$\ov{\phi}\colon N\to N$ and a birational morphism $\epsilon\colon
N/\langle\ov{\phi}\rangle\to W$ with finite fibers such that $\ov{f}$
is the composition
\begin{equation}
N\overset{\pi}{\lra} N/\langle\ov{\phi}\rangle
\overset{\epsilon}{\lra} W.
\end{equation}
In particular $N$ is projective and $\deg W=6$.
\end{lmm}
\begin{proof}
Let's show that $\ov{f}$ has finite fibers of cardinality at most $2$ and that the generic fiber has cardinality $2$.  The fiber of $\ov{f}$ over
$|I_{F^3_5}(2)|\in |I_S(2)|^{\vee}$ consists of the unique singular point $p$ of $N$. Now let $\Lambda\in(W\setminus\{|I_{F^3_5}(2)|\})$, i.e.
\begin{equation}\label{parcondicio}
\Lambda=f([Z]),\quad [Z]\notin P.
\end{equation}
Then $\ov{f}^{-1}(\Lambda)=f^{-1}(\Lambda)$.
Let $\Lambda_0:=\Lambda\cap |I_{F^3_5}(2)|$ and choose
$\lambda_0\in(\Lambda\setminus |I_{F^3_5}(2)|)$.
One has
\begin{equation}
\bigcap\limits_{\lambda\in\Lambda_0} Q_{\lambda}=
 F^3_5\cup A_{\Lambda}
\end{equation}
where $A_{\Lambda}$ is a plane such that $A_{\Lambda}\cap
F^3_5=C_{\Lambda}$ is a conic. (See the proof of Lemma~(4.20)
of~\cite{og1}.) We claim that $Q_{\lambda_0}\not\supset A_{\Lambda}$:
in fact if $Q_{\lambda_0}\supset A_{\Lambda}$ then $C_{\Lambda}\subset
S$ contradicting~(\ref{retteconiche}). Thus  $Q_{\lambda_0}\cap
A_{\Lambda}$ is a conic $C'_{\Lambda}$. By~(\ref{parcondicio})  the
line $span(Z)$ is contained in $C'_{\Lambda}$. Thus $C'_{\Lambda}$ is
degenerate and $f^{-1}(\Lambda)$ consists of the set of lines contained
in $C'_{\Lambda}$. This shows that
$f^{-1}(\Lambda)=\ov{f}^{-1}(\Lambda)$ has cardinality at most $2$. It
also follows easily that the generic fiber of $\ov{f}$ consists of $2$
distinct points. Since $N$ is normal there is a regular covering
involution $\ov{\phi}$ such that $\ov{f}$ factors through the quotient
map $N\to N/\langle\ov{\phi}\rangle$. Since $\ov{f}$ has finite fibers
so does $\epsilon$, since the generic fiber of $\ov{f}$ consists of $2$
points the map $\epsilon$ is birational. The  line-bundle
$\ov{f}^{*}\cO_W(1)$ is ample because $\ov{f}\colon N\to W$ has finite
fibers, thus  $N$ is projective. We know (see~(4.3) of~\cite{og1}) that
\begin{equation}
\int_{S^{[2]}}c_1(f^{*}\cO_W(1))=12.
\end{equation}
Since $f\colon S^{[2]}\to W$ has generic fiber of cardinality $2$ we
get that $\deg W=6$.
\end{proof}
We will show that the map $\epsilon$ of Lemma~(\ref{gradoeffe}) is an
isomorphism.  Since $(S^{[2]}\setminus P)\cong (N\setminus\{p\})$ the
involution $\ov{\phi}$  defines a birational involution
\begin{equation}\label{mappafi}
\phi\colon S^{[2]}\dashrightarrow S^{[2]}.
\end{equation}
(This is the birational involution of Proposition~(4.21) of~\cite{og1}.)
 Let
\begin{equation}
 \beta\colon Bl_P(S^{[2]})\to S^{[2]}
\end{equation}
 be the blow-up of $P$:  since $P$
 is Lagrangian the symplectic form on $S^{[2]}$ induces an isomorphism
$N_{P/S^{[2]}}\cong\Omega^1_P$ and hence the exceptional divisor of
$\beta$ is identified with
 the incidence
variety $\Gamma\subset P\times P^{\vee}$  and  the restriction of
$\beta$ to the exceptional divisor is identified with the projection
$\Gamma\to P$. We abuse notation and view $\Gamma$ as the exceptional
divisor in $Bl_P(S^{[2]})$
\begin{lmm}\label{sollevofi}
The map $\phi$ of~(\ref{mappafi}) is not regular along $P$. There is a
regular involution $\wt{\phi}\colon Bl_P(S^{[2]})\to Bl_P(S^{[2]})$
which is equal to $\phi$ on $(S^{[2]}\setminus P)\subset
Bl_P(S^{[2]})$. There is an identification $P\cong P^{\vee}$ such that
$\wt{\phi}|_{\Gamma}$ is induced by the involution on $P\times
P^{\vee}$ which interchanges the factors.
\end{lmm}
\begin{proof}
The eigenspaces of  the isometry $H^2(\phi)$ on $H^2(S^{[2]})$ induced
by $\phi$ are given by (see~(4.3) of~\cite{og1})
\begin{equation}\label{azione}
 H^2(\phi)_{+}=\CC f^{*}\cO_W(1),\quad
 H^2(\phi)_{-}=f^{*}\cO_W(1)^{\bot}.
\end{equation}
Suppose $\phi$ is regular along $P$: if $D$ is an ample divisor on
$S^{[2]}$ then $c_1(D+\phi^{*}D)$ is an ample $\phi$-invariant class,
this contradicts~(\ref{azione}) because $f^{*}\cO_W(1)$ is not ample.
Let $\psi\colon S^{[2]}\dashrightarrow \cX$ be the flop of $P$: thus
$\psi$ is the inverse of the blowup $Bl_P(S^{[2]})\to S^{[2]}$ followed
by the morphism $Bl_P(S^{[2]})\to\cX$ which contracts $\Gamma$ along
the \lq\lq other\rq\rq fibration $\Gamma\to P^{\vee}$. In particular
$\cX$ contains $P^{\vee}$. Let $\ell\subset P$ and $\ell^{\vee}\subset
P^{\vee}$ be lines. The isometry $H^2(\psi)$ identifies
$(\ell^{\vee})^{\bot}$ with $\ell^{\bot}$. In fact we have contractions
$c\colon S^{[2]}\to N$ and $c^{\vee}\colon\cX\to N$ which give
identifications $\ell^{\bot}=H^2(N)=(\ell^{\vee})^{\bot}$. In
particular $H^2(\psi)$ sends a nef divisor class in
$(\ell^{\vee})^{\bot}$ to a nef divisor class. On the other hand
$H^2(\psi)$ maps the half-space $\ell^{\vee}_{>0}$ to the half-space
$\ell_{<0}$. By~(\ref{azione}) we get that $(\psi\circ\phi)^{*}$ maps
an ample divisor to an ample divisor: since $(\psi\circ\phi)$ defines a
regular map between the complements of subsets of codimension $2$  it
follows that $(\psi\circ\phi)$ is regular and hence an isomorphism. It
follows also that $\phi$ induces a regular involution $\wt{\phi}\colon
Bl_P(S^{[2]})\to Bl_P(S^{[2]})$. Let's show that the restriction of
$\wt{\phi}$ to $\Gamma$ is as stated. Any automorphism of $\Gamma$
sends the projection $\Gamma\to P$ to itself composed with an
automorphism of $P$ or to the projection $\Gamma\to P^{\vee}$ composed
with an automorphism of $P^{\vee}$. Since $\phi$ is not regular the
latter holds and it follows that $\wt{\phi}|_{\Gamma}$ is as stated.
\end{proof}
\begin{crl}\label{proprioquoz}
The map $\epsilon\colon N/\langle\ov{\phi}\rangle\to W$ of
Lemma~(\ref{gradoeffe}) is an isomorphism.
\end{crl}
\begin{proof}
The quotient $Bl_P(S^{[2]})/\langle\wt{\phi}\rangle$ is a projective
birational model of $N/\langle\ov{\phi}\rangle$ and hence it is
birational to $W$ by Lemma~(\ref{gradoeffe}). The Kodaira dimension of
$Bl_P(S^{[2]})/\langle\wt{\phi}\rangle$ is $0$ hence also the Kodaira
dimension of $W$ is $0$. By Lemma~(\ref{gradoeffe}) we know that $\deg
W=6$ and hence by adjunction $W$ is smooth in codimension $1$. Thus $W$
is normal: since $\epsilon$ is birational with finite fibers we get
that $\epsilon$ is an isomorphism.
\end{proof}
Now we reintroduce the subscript $t$; thus we have $S_t$, $f_t$, $N_t$,
$W_t$ etc.
\begin{prp}\label{yew}
Keep notation as above. Let $A\in\Delta^0(V)$ and let $t\in T'''$ such
that $X_{A_t}\cong M_t$. Thus by Proposition~(\ref{emmetiepw}) and
Corollary~(\ref{nontutto}) we have
$Y_{A^{\bot}}\subset|I_{S_t}(2)|^{\vee}$. Then
\begin{equation}
Y_{A^{\bot}}= W_t.
\end{equation}
\end{prp}
\begin{proof}
Proposition~(5.20) of~\cite{og1} gives that the reduced scheme
$(Y_{A^{\bot}})_{red}$ is equal to $W_t$ (the hypothesis of that
Proposition is that~(\ref{dieci}) holds, however the same proof goes
through because all that is needed is the validity
of~(\ref{ipsilonti})). Now $Y_{A^{\bot}}$ is a degree-$6$ hypersurface
because $Y_{A^{\bot}}\not=\PP(V^{\vee})$ and on the other hand $W_t$ is
a degree-$6$ hypersurface by Lemma~(\ref{gradoeffe}) and hence from
$(Y_{A^{\bot}})_{red}=W_t$ we get that $Y_{A^{\bot}}$ is reduced and
equal to $W_t$.
\end{proof}
Let $A\in\Delta^0(V)$ and let $t\in T'''$ such that $X_{A_t}\cong M_t$.
Then $|I_{F^3_5}(2)|\in W_t$ - see~(\ref{puntospeciale} - and hence by
the above proposition $|I_{F^3_5}(2)|\in Y_{A^{\bot}}$; we denote this
point by $q_{A^{\bot}}$.
\begin{prp}\label{descrivoduale}
Let $A\in\Delta^0(V)$ and let $t\in T'''$ such that $X_{A_t}\cong M_t$.
Then
\begin{equation}\label{molteplicita}
{\rm mult}_{q_{A^{\bot}}}Y_{A^{\bot}}=3.
\end{equation}
Let $y\in (Y_{A^{\bot}}\setminus\{q_{A^{\bot}}\})$ and hence
$\ov{f_t}^{-1}(y)=f_t^{-1}(y)$ consists of two points or of one point.
\begin{itemize}
\item [(a)]
If $f_t^{-1}(y)$ consists of two points then $Y_{A^{\bot}}$ is smooth
at $y$.
\item [(b)]
If $f_t^{-1}(y)$ consists of a single point then the analytic germ
$(Y_{A^{\bot}},y)$  is isomorphic to the product of a smooth
$2$-dimensional germ and the germ of an $A_1$-singularity.
\end{itemize}
\end{prp}
\begin{proof}
By Proposition~(6.2) of~\cite{og2}  $Y_A^{\vee}$ has multiplicity $3$
in $|I_{F^3_5}(2)|$; this proves (\ref{molteplicita}). In order to
prove~(a)-(b) we notice that  $Y_{A^{\bot}}=W_t$ and by
Proposition~(\ref{proprioquoz}) the map $\ov{f}_t\colon N_t\to W_t$ is
identified with the quotient map $N_t\to N_t/\langle\ov{\phi}\rangle$.
Since $p_t=\ov{f}^{-1}(q_A)$ is the unique singular point of $N_t$ this
gives Item~(a). In order to prove Item~(b)  we notice that
$(N_t\setminus\{p_t\})=(S_t^{[2]}\setminus P_t)$ and that the quotient
map $(N_t\setminus\{p_t\})\to
(N_t\setminus\{p_t\})/\langle\ov{\phi}\rangle$ is identified with the
quotient map $(S_t^{[2]}\setminus P_t)\to (S_t^{[2]}\setminus
P_t)/\langle\phi\rangle$. By Proposition~(4.21) of~\cite{og1} the
restriction of $\phi$ to $(S_t^{[2]}\setminus P_t)$ is an
anti-symplectic involution and hence its fixed point set is a
Lagrangian surface; this proves~(b).
\end{proof}
Let $A\in\Delta^0(V)$; then by (\ref{molteplicita})
 the sextic $Y_{A^{\bot}}$ has a point of multiplicity $3$ and hence
 $A^{\bot}\notin\mathbb{LG}(\wedge^3 V^{\vee})^0$ because if
 $B\in\mathbb{LG}(\wedge^3 V^{\vee})^0$ then $Y_B$ has points of
 multiplicity at most $2$. Thus
 \begin{equation}
    \Delta^0(V)\cap\mathbb{LG}(\wedge^3 V)^{00}=\es.
\end{equation}
Let $\Delta^0_{*}(V)\subset\Delta^0(V)$
 be the set of  points which are smooth points of
the projective variety $\left(\mathbb{LG}(\wedge^3
V)\setminus\mathbb{LG}(\wedge^3 V)^{00}\right)$. By
Proposition~(\ref{codelta})
 we know that $\Delta^0(V)$ is locally (in the classical topology) a
 codimension $1$ submanifold of $\mathbb{LG}(\wedge^3 V)$;
 since $\Delta^0(V)\cap \mathbb{LG}(\wedge^3 V)^{00}=\es$ it
 follows that  $\Delta^0_{*}(V)$ is open dense in
 $\Delta^0(V)$. We let $T_{*}\subset T'''$ be the set of $t$ such that $M_t\cong X_A$ for some $A\in \Delta^0_{*}(V)$; since $\Delta^0_{*}(V)$ is open dense in
 $\Delta^0(V)$ also $T_{*}$ is open dense in $T'''$.
 Let
\begin{equation}
\mathbb{LG}(\wedge^3 V)^0_{*}:=\mathbb{LG}(\wedge^3
V)^{00}\cup\Delta^0_{*}(V).
\end{equation}
\begin{prp}
Keep notation as above. Then $\mathbb{LG}(\wedge^3 V)^0_{*}$ is open in
$\mathbb{LG}(\wedge^3 V)$ (for the classical topology) and
$\Delta^0_{*}(V)$ is a non-empty codimension-$1$ submanifold of
$\mathbb{LG}(\wedge^3 V)^0_{*}$.
\end{prp}
\begin{proof}
Let $A\in \Delta^0_{*}(V)$. By definition of $\Delta^0_{*}(V)$ there
exists an open $U\subset\mathbb{LG}(\wedge^3 V)$ such that
$U\cap\Delta^0_{*}(V)=U\cap \left(\mathbb{LG}(\wedge^3
V)\setminus\mathbb{LG}(\wedge^3 V)^{00}\right)$ and hence $U\subset
\mathbb{LG}(\wedge^3 V)^0_{*}$; this proves that $\mathbb{LG}(\wedge^3
V)^0_{*}$ is open. We have already proved that $\Delta^0_{*}(V)$ is a
non-empty locally closed codimension-$1$ submanifold  of
$\mathbb{LG}(\wedge^3 V)$; this gives the second statement of the
proposition.
\end{proof}
We let
\begin{align}
\Delta^\infty_{*}(V^{\vee})& :=\delta(\Delta^0_{*}(V)), \\
\mathbb{LG}(\wedge^3 V^{\vee})^\infty_{*}& :=\delta(
 \mathbb{LG}(\wedge^3 V)^0_{*}). \\
\end{align}
Of course every definition above has a \lq\lq dual \rq\rq definition
obtained by substituting $V^{\vee}$ to $V$; thus we have
$\Delta^0_{*}(V^{\vee})\subset \mathbb{LG}(\wedge^3 V^{\vee})^0$,
$\Delta^\infty_{*}(V)\subset \mathbb{LG}(\wedge^3 V)$ etc. Let
\begin{multline}
\mathbb{LG}(\wedge^3 V^{\vee})^\sharp  :=
 \mathbb{LG}(\wedge^3 V^{\vee})^0_{*}\cup
 \mathbb{LG}(\wedge^3 V^{\vee})^\infty_{*}= \\
 =\mathbb{LG}(\wedge^3 V^{\vee})^{00}\cup\Delta^0_{*}(V^{\vee})
 \cup\Delta^\infty_{*}(V^{\vee}).
\end{multline}
 Let
$\cY(V^{\vee})\subset\mathbb{LG}(\wedge^3
V^{\vee})^\sharp\times\PP(V^{\vee})$ be the tautological EPW-sextic;
thus $\cY(V^{\vee})\cap\{B\}\times\PP(V^{\vee})=Y_B$. If
$B\in\Delta^\infty_{*}(V^{\vee})$ then by Proposition
(\ref{descrivoduale}) there is a unique point $q_B\in Y_B$ of
multiplicity strictly greater than $2$; let
\begin{align}
\cQ(V^{\vee}) & :=  \{(B,q_B)|\ B\in\Delta^\infty_{*}(V^{\vee})\}, \\
\cY(V^{\vee})^{\sharp} & := \cY(V^{\vee})\setminus\cQ(V^{\vee}).
\end{align}
\begin{prp}\label{doppio}
There exists a double cover
$f\colon\cX(V^{\vee})\to\cY(V^{\vee})^{\sharp}$ with the following
properties.
\begin{itemize}
\item [(1)]
Let $\pi\colon \cX(V^{\vee})\to\mathbb{LG}(\wedge^3 V^{\vee})^\sharp$  be the composition of $f$ and the projection
$\cY(V^{\vee})^{\sharp}\to\mathbb{LG}(\wedge^3 V^{\vee})^\sharp$; then $\pi$ is a submersion of smooth manifolds.
\item [(2)]
Let $B\in\mathbb{LG}(\wedge^3 V^{\vee})^0_{*}$; then $\pi^{-1}(B)\cong
X_B$ and the map $\pi^{-1}(B)\to Y_B$ defined by $f$ is isomorphic to
the natural double cover $X_B\to Y_B$.
\item [(3)]
Let $B\in\Delta^\infty_{*}(V^{\vee})$ and $A:=B^{\bot}$. Since
$A\in\Delta^0_{*}(V)$ there exists $t\in T_{*}$ such that $X_A\cong M_t$
and $Y_A\cong\Sigma'_t$. Then $\pi^{-1}(B)\cong
(S^{[2]}_t\setminus\{P_t\})$ and the map $\pi^{-1}(B)\to
(Y_B\setminus\{q_B\})$ defined by $f$ is isomorphic to the double cover
$(S^{[2]}_t\setminus\{P_t\})\to (W_t\setminus |I_{F^3_5}(2)|)$ given by
the restriction of~(\ref{mappaeffe}). (Recall that
$(Y_B\setminus\{q_B\})\cong (W_t\setminus |I_{F^3_5}(2)|)$ by
Proposition~(\ref{yew})).
\end{itemize}
\end{prp}
\begin{proof}
 Let $B\in\mathbb{LG}(\wedge^3 V^{\vee})^{\sharp}$. Then
 $Y_B\not=\PP(V^{\vee})$ by Corollary~(\ref{nontutto}); thus
 the map $\lambda_B$ defined in Section~(\ref{prologo}) (with $V^{\vee}$ replacing $V$) is non-zero  and $Y_B$ is the zero-scheme of $\det(\lambda_B)$. Since $Y_B$ is a Lagrangian degeneracy locus there exists
 locally in $\PP(V^{\vee})$  a symmetric map of
 vector-bundles  giving a resolution of
$coker(\lambda_B)$, i.e.~we can cover $\PP(V^{\vee})$ by open sets $U$ such that on each $U$ we have a locally-free resolution
\begin{equation}
0\to E_U\overset{\alpha_U}{\lra} E_U^{\vee}\lra
 coker(\lambda_B)|_U\to 0
\end{equation}
where $\alpha_U$ is a symmetric map of vector-bundles.
Furthermore if $F\hra\wedge^3 V^{\vee}\otimes\cO_{\PP(V^{\vee})}$ is the Lagrangian sub-vector-bundle defined in Section~(\ref{prologo})  (with $V^{\vee}$ replacing $V$) then
there is an isomorphism $F_p\cap B\cong ker(\alpha_p)$ for all $p\in U$. Thus if $\dim(F_p\cap B)\ge r$ then $mult_p(Y_B)\ge r$. Let $p\not=q_B$. Then  $mult_p(Y_B)\le 2$  by Proposition~(\ref{descrivoduale}) and hence $\dim(F_p\cap B)\le 2$. Furthermore one of the following holds:
\begin{itemize}
\item[(1)]
If  $\dim(F_p\cap B)=1$ then locally around $p$ we have $coker(\lambda_B)\cong i_{*}\cO_{Y_B}$ where $i\colon Y_B\hra\PP(V^{\vee})$ is the inclusion.
\item[(2)]
If  $\dim(F_p\cap B)=2$ there exist an open (in
the classical topology) $U\subset \PP(V^{\vee})$  containing $p$, functions
 $x,y,z\in Hol(U)$ vanishing at $p$ with
  linearly independent differentials and an
 exact sequence 
\begin{equation}\label{risolvo}
  0\to \cO^2_U\mapor{M}\cO^2_U\mapor{}
  coker(\lambda_B)|_U\to 0
\end{equation}
  where $M$ is the map defined by the matrix
\begin{equation}
\left(\begin{array}{rr}
     x &  y \\
    y & z
  \end{array}\right).
\end{equation}
\end{itemize}
In particular we see that there exists a sheaf $\zeta_B$ on $(Y_B\setminus\{q_B\})$ such that
outside $q_B$ we have $coker(\lambda_B)=i^0_{*}\zeta_B$ where   $i^0\colon (Y_B\setminus\{q_B\})\hra(\PP(V^{\vee})\setminus\{q_B\})$ is the inclusion. From the local description of  $coker(\lambda_B)$ given above we also get a canonical isomorphism of sheaves on $(\PP(V^{\vee})\setminus\{q_B\})$:
\begin{equation}\label{dualeduale}
Ext^1(coker(\lambda_B),\cO_{\PP(V^{\vee})})
 |_{\PP(V^{\vee})\setminus\{q_B\}}
 \cong i^0_{*}(\zeta_B^{\vee}\otimes N_{Y_B/\PP(V^{\vee})}).
\end{equation}
 (See Proposition~(4.3) of~\cite{og2}.) 
Let  $B^{\vee}\subset\wedge^3 V^{\vee}$  be a Lagrangian subspace complementary to $B$; thus we have a direct sum decomposition 
\begin{equation}\label{decovu}
\wedge^3 V=B\oplus B^{\vee}. 
\end{equation}
Then $\lambda_B$ can be identified with the map of vector-bundles $F\to B^{\vee}\otimes\cO_{\PP(V^{\vee})}$ associated to Decomposition~(\ref{decovu}). Let
$\mu_B\colon F\to B\otimes\cO_{\PP(V^{\vee})}$ be the \lq\lq other\rq\rq map associated to Decomposition~(\ref{decovu}). The diagram
\begin{equation}
 \begin{array}{ccc}
F&\mapor{\lambda_B}&
B^{\vee}\otimes\cO_{\PP(V^{\vee})}\\
\mapver{\mu_B}& &\mapver{\mu_B^{\vee}}\\
B\otimes\cO_{\PP(V^{\vee})}& \mapor{\lambda_B^{\vee}}&
F^{\vee}
\end{array}
\end{equation}
is commutative because
$F\overset{(\mu_B,\lambda_B)}{\lra}(B\oplus
B^{\vee})\otimes\cO_{\PP(V^{\vee})}$ is a Lagrangian
embedding. The map $\lambda_B$ is an injection
of sheaves because $Y_B\not=\PP(V^{\vee})$ and hence
also $\lambda_B^{\vee}$ is an injection of
sheaves. Thus there is a unique $\beta_B\colon
coker(\lambda_B)\lra
Ext^1(coker(\lambda_B),\cO_{\PP(V^{\vee})})$ making the
following diagram commutative with exact
horizontal sequences:
\begin{equation}\label{spqr}
\begin{array}{ccccccccc}
0 & \to & F&\mapor{\lambda_B}& B^{\vee}\otimes\cO_{\PP(V^{\vee})} & \lra & coker(\lambda_B)
&
\to & 0\\
 & & \mapver{\mu_B}& &\mapver{\mu^{\vee}_B} &
&
\mapver{\beta_B}& & \\
0 & \to & B\otimes\cO_{\PP(V^{\vee})}& 
\mapor{\lambda_B^{\vee}}& F^{\vee} & \lra &
Ext^1(coker(\lambda_B),\cO_{\PP(V^{\vee})}) & \to & 0
\end{array}
\end{equation}
By Isomorphism~(\ref{dualeduale}) we get that the restriction of 
$\beta_B$ to $(\PP(V^{\vee})\setminus\{q_B\})$ defines a map of sheaves on $(Y_B\setminus\{q_B\})$
\begin{equation}\label{zetazeta}
\zeta_B\to \zeta_B^{\vee}(6).
\end{equation}
Since $F\overset{(\mu_B,\lambda_B)}{\lra}(B\oplus
B^{\vee})\otimes\cO_{\PP(V^{\vee})}$ is an injection of vector-bundles the above map is an isomorphism - this follows from Claim~(4.5) of~\cite{og2}.
Let $\xi_B:=\zeta_B(-3)$; then~(\ref{zetazeta}) defines an isomorphism $\xi_B\overset{\sim}{\lra}\xi_B^{\vee}$ which is symmetric and hence it gives $\cO_{(Y_B\setminus\{q_B\})}\oplus\xi_B$ the structure of a commutative finite $\cO_{(Y_B\setminus\{q_B\})}$-algebra. 
Let $X_B:= Spec(\cO_{(Y_B\setminus\{q_B\})}\oplus\xi_B)$ and 
$f_B\colon X_B\to 
(Y_B\setminus\{q_B\})$ be the structure map: clearly $f_B$ is finite of degree $2$. The above construction is the analogue of the construction of the natural double cover $X_B\to Y_B$ for $B\in\mathbb{LG}(\wedge^3 V^{\vee})^0$; thus we have
 a double cover
$f\colon\cX(V^{\vee})\to\cY^{\sharp}(V^{\vee})$ such that Item~(2) holds and such that for $B\in\Delta^{\infty}_{*}(V^{\vee})$ we have $\pi^{-1}(B)=X_B$ and the map $\pi^{-1}(B)\to (Y_B\setminus\{q_B\})$ 
is the structure map $f_B$ defined above. It remains to prove that Items~(1) and~(3) hold. Let
\begin{align}
C_B:= & \{p\in(sing Y_B\setminus q_B)|\ \dim(F_p\cap B)=1\} \\
D_B:= & \{p\in(sing Y_B\setminus q_B)|\ \dim(F_p\cap B)=2\}.
\end{align}
Thus $(sing Y_B\setminus q_B)=C_B\coprod D_B$. 
The local description given above of $coker(\lambda_B)$ near $p\in(Y_B\setminus q_B)$ shows that both $C_B$ and $D_B$ are smooth and closed (in $(Y_B\setminus q_B)$). Let $p\in (Y_B\setminus q_B)$; the map $f_B$ behaves differently depending on whether $p\in C_B$ or $p\in D_B$: in fact (see~\cite{og2})
\begin{itemize}
\item[($\alpha$)]
If $p\in C_B$ then $f_B$ is unramified over $p$.
 \item[($\beta$)]
If $p\in D_B$ then $f_B$ is ramified over $p$ and $X_B$ is smooth at $f_B^{-1}(p)$.
\end{itemize}
Let's prove that 
\begin{equation}\label{zorro}
X_B\cong(S_t^{[2]}\setminus P_t).
\end{equation}
Let $B'\in \mathbb{LG}(\wedge^3 V^{\vee})^0_{*}$ and define $D_{B'}$ as above. Then $D_{B'}$ is not empty in fact it is a surface; it follows that $D_B\not=\es$ and by Item~($\beta$) above we get that the \'etale covering
\begin{equation}\label{valentino}
\left(X_B\setminus f_B^{-1}(sing Y_B)\right)\lra 
 (Y_B\setminus sing Y_B)
\end{equation}
is not trivial. On the other hand let  $t\in T_{*}$ such that 
$X_{A^{\bot}}\cong M_t$
and $Y_{A^{\bot}}\cong\Sigma'_t$. 
Let $f_t\colon S_t^{[2]}\to Y_B$
be the map defined by~(\ref{mappaeffe}); then $f_t$ defines an \'etale double cover
\begin{equation}\label{agostini}
(S_t^{[2]}\setminus f_t^{-1}(sing Y_B))\lra 
 (Y_B\setminus sing Y_B).
\end{equation}
Now $S_t^{[2]}$ is  simply connected and $f_t^{-1}(sing Y_B)$ has codimension $2$ in $S_t^{[2]}$ hence $(S_t^{[2]}\setminus f_t^{-1}(sing Y_B))$ is simply connected. Thus $\pi_1(Y_B\setminus sing Y_B)\cong\ZZ/(2)$ and hence there is a unique non-trivial double cover of $(Y_B\setminus sing Y_B)$ and it is given by~(\ref{agostini}). Since~(\ref{valentino}) is a non-trivial double cover it follows that $C_B=\es$ and that~(\ref{zorro}) holds. This proves Items~(1)-(3) of the proposition.
\end{proof}
\subsection{Extension of the local period map across $\Delta^{\infty}_{*}(V^{\vee})$}
\label{verifica}
\setcounter{equation}{0}
We will prove that the local  period map
extends across $\Delta^{\infty}_{*}(V^{\vee})$.  
Let $\pi^0\colon\cX(V^{\vee})^0_{*}\to\mathbb{LG}(\wedge^3 V^{\vee})^0_{*}$ be the tautological family of double EPW-sextics i.e.
the restriction of the map $\pi$ of Proposition~(\ref{doppio}) to $\pi^{-1}(\mathbb{LG}(\wedge^3 V^{\vee})^0_{*})$. Since $\pi^0$ is proper it defines the variation of Hodge structures
\begin{equation}\label{variazione}
(R^2\pi^0_{*}\ZZ,F^p)
\end{equation}
where $F^0\supset F^1\supset F^2$ is the Hodge filtration of $(R^2\pi^0_{*}\ZZ)\otimes\cO_{\mathbb{LG}(\wedge^3 V^{\vee})^0_{*}}$ i.e.~the fiber of $F^p$ over $A$ is 
\begin{equation}
F^p_A=F^p H^2(X_A):=\bigoplus_{p'\ge p} H^{p',2-p'}(X_A).
\end{equation}
 Furthermore we have the symmetric section 
 \begin{equation}\label{sezionebibi}
\BB\in H^0((R^2\pi^0_{*}\QQ)\otimes(R^2\pi^0_{*}\QQ))
\end{equation}
which gives the Beauville-Bogomolov bilinear form on $H^2(X_A)$ for every $A\in \mathbb{LG}(\wedge^3 V^{\vee})^0_{*}$. 
 \begin{prp}\label{estensione}
There exist
\begin{itemize}
\item[(1)]
a local system $\HH(V^{\vee})$ on $\mathbb{LG}(\wedge^3 V^{\vee})^{\sharp}$,
\item[(2)]
 a decreasing filtration of $\HH(V^{\vee})\otimes\cO_{\mathbb{LG}(\wedge^3 V^{\vee})^{\sharp}}$ by holomorphic sub-bundles $\wt{F}^0\supset \wt{F}^1\supset \wt{F}^2$,
\item[(3)]
a symmetric section $\wt{\BB}\in H^0((\HH(V^{\vee})\otimes\HH(V^{\vee}))\otimes\QQ))$
\end{itemize}
such that the following hold:
\begin{itemize}
\item[(a)]
$(\HH(V^{\vee}),\wt{F}^p)$ extends the variation of Hodge structures~(\ref{variazione}). 
\item[(b)]
 Let $B\in\Delta^{\infty}_{*}(V^{\vee})$ and $t\in T_{*}$ such that $M_t\cong X_{B^{\bot}}$; then there exists an isomorphism of Hodge structures  
 \begin{equation}\label{forzaitalia}
 (\HH(V^{\vee})_B,\wt{F}^p_B)\cong H^2(S_t^{[2]}).
\end{equation}
\item[(c)]
The restriction of $\wt{\BB}$ to $\mathbb{LG}(\wedge^3 V^{\vee})^0_{*}$ is equal to $\BB$.
\item[(d)]
 Let $B\in\Delta^{\infty}_{*}(V^{\vee})$ and $t\in T_{*}$ such that $M_t\cong X_{B^{\bot}}$; then  Isomorphism~(\ref{forzaitalia}) is an isometry between  $(\HH(V^{\vee})_B,\wt{\BB}_B)$ and 
 $H^2(S_t^{[2]})$ equipped with the Beauville-Bogomolov symmetric bilinear form.
\end{itemize}
\end{prp}
\begin{proof}
Let 
$\cU\subset\mathbb{LG}(\wedge^3 V^{\vee})^{\sharp}$ be an  open ball. We assume that $\cU$ is small: then there exists a hyperplane $H\subset\PP(V^{\vee})$
such that $q_B\notin H$ for all $B\in
\cU\cap\Delta^{\infty}_{*}(V^{\vee})$. Let $\pi$ and $f$ be as in Proposition~(\ref{doppio}); we let
\begin{align}
\cX(V^{\vee})_{\cU}:= & \pi^{-1}(\cU), \\
 Z({\cU},H):= & \cX(V^{\vee})_{\cU}\cap f^{-1}H.
\end{align}
Let $\rho\colon Z({\cU},H)\to\cU$ be given by the restriction of $f$. By our choice of $H$ the map $\rho$ is  proper submersive with fibers smooth $3$-folds. Thus we have a variation of Hodge structures $(R^2\rho_{*}\ZZ,F^p)$; we denote it by $(\HH(\cU,H),F^p(\cU,H))$. Let $B\in(\cU\setminus\Delta^{\infty}_{*})$; then $rho^{-1}(B)$ is an ample divisor on $X_B$ and hence we have a canonical isomorphism
\begin{equation}
(\HH(\cU,H),F^p(\cU,H))|_{(\cU\setminus\Delta^{\infty}_{*})}
\cong (R^2\rho_{*}\ZZ,F^p)|_{(\cU\setminus\Delta^{\infty}_{*})}.
\end{equation}
 This shows that $(\HH(\cU,H),F^p(\cU,H))$ does not depend on the choice of $H$ and that the collection of $(\HH(\cU,H),F^p(\cU,H))$ gives an extension $(\HH(V^{\vee}),\wt{F}^p)$ of the variation of Hodge structures~(\ref{variazione}). Now let's prove Item~(b). Let $B\in(\cU\cap\Delta^{\infty}_{*}(V^{\vee}))$ and let $t$ be as in Item~(b). Let $f_t\colon S_t^{[2]}\to W_t=Y_B$ be the map given by~(\ref{mappaeffe}). Then $\rho^{-1}(B)=f_t^{-1}H$ and since $f_t$ is semi-small 
the restriction map 
$H^2(S_t^{[2]})\to H^2(f_t^{-1}H)$ is
an isomorphism of (integral) Hodge structures; this proves Item~(b) because $(\HH(V^{\vee})_B,\wt{F}^p)$ is isomorphic to the Hodge structure on $H^2(f_t^{-1}H)$ by definition. We define $\wt{\BB}$ as follows. Let $\cU$ and $\rho\colon Z({\cU},H)\to\cU$ be as above; let $Z_B:=\rho^{-1}(B)$. For $B\in\cU$ we have the Lefschetz decomposition
\begin{equation}\label{lefdec}
\HH(V^{\vee})_B= H^2(Z_B;\QQ)=\QQ c_1(L_B)|_{Z_B}\oplus 
H^2(Z_B;\QQ)_{prim}
\end{equation}
where $L_B$ is the tautological ample line-bundle on $X_B$. 
We let $\wt{\BB}(\cU,H)_B$  be the symmetric bilinear form on $\HH(V^{\vee})_B$ characterized by the following requirements: 
\begin{itemize}
\item[($\alpha$)]
Decomposition~(\ref{lefdec}) is orthogonal for $\wt{\BB}(\cU,H)_B$.
\item[($\beta$)]
$\wt{\BB}_B(c_1(L_B)|_{Z_B},c_1(L_B)|_{Z_B})=2$.
\item[($\gamma$)]
 If $x,y\in H^2(Z_B;\QQ)_{prim}$ then $\wt{\BB}(\cU,H)_B (x,y)=
 \frac{1}{2}\int_{Z_B}c_1(L_B)\wedge x\wedge y$.
\end{itemize}
Since the Lefscethz decomposition is flat for the Gauss-Manin connection we have a well-defined section $\wt{\BB}(\cU,H)\in H^0(\HH(V^{\vee})|_{\cU})$ with value $\wt{\BB}(\cU,H)_B$ at $B\in \cU$. 
Let $B\in(\cU\setminus\Delta^{\infty}_{*}(V^{\vee})$ and let $\iota\colon Z_B\hra X_B$ be the inclusion. If $\xi_1,\xi_2,\xi_3\in H^2(X_B)$ then
\begin{equation}\label{polarizzo}
\int_{Z_B}\iota^{*}\xi_1\wedge\iota^{*}\xi_2\wedge
\iota^{*}\xi_3=
\int_{X_B}c_1(L_B)\wedge\xi_1\wedge\xi_2\wedge\xi_3= 
\frac{1}{2}\sum_{\sigma\in\cS_3}(c_1(L_B),\xi_{\sigma(1)})\cdot (\xi_{\sigma(2)},\xi_{\sigma(3)}).
\end{equation}
It follows from this that $\wt{\BB}(\cU,H)$ does not depend on $H$ and that the collection of $\wt{\BB}(\cU,H)$'s defines an extension of $\BB$. Item~(d) holds because Formula~(\ref{polarizzo}) holds if we replace $X_B$ by $S_t^{[2]}$ . 
\end{proof}
The map $\delta$ defines an isomorphism 
\begin{equation}
\begin{matrix}
\mathbb{LG}(\wedge^3 V)^{\sharp} & \overset{\sim}{\lra} &
\mathbb{LG}(\wedge^3 V)^{\sharp} \\
A & \mapsto & \delta(A)=A^{\bot}.
\end{matrix}
\end{equation}
From now on we will denote by $\delta$ what is actually the restriction of $\delta$ to $\mathbb{LG}(\wedge^3 V)^{\sharp}$.
Let  $\cU\subset\mathbb{LG}(\wedge^3 V^{\vee})^{\sharp}$ be an open set: then we have  two local systems, namely 
\begin{equation}\label{sistloc}
\HH(V)|_{\cU},\qquad \HH(V^{\vee})|_{\delta(\cU)}.
\end{equation}
\begin{dfn}\label{marchio}
Assume that both local systems~(\ref{sistloc}) are trivial.
A {\it marking of $\HH(V)|_{\cU}$} is an isomorphism $\Psi\colon\HH(V)|_{\cU}\overset{\sim}{\lra}\cU\times\wt{\Lambda}$ such that:
\begin{itemize}
\item[(1)]
 $\wt{\BB}_{A}(x,y)=(\Psi(x),\Psi(y))$ for every $A
\in\cU$ and $x,y\in \HH(V)_{A}$, and
\item[(2)]
$\Psi$ sends the flat section $A\mapsto  c_1(L_{A})|_{Z(\cU,H)}$ (notation as in the proof of Proposition~(\ref{estensione})) to $\cU\times u$ where $u$ is given by~(\ref{scelgou}).
\end{itemize}
We define similarly a {\it marking $\Psi$ of $\HH(V^{\vee})|_{\delta(\cU)}$}.  
\end{dfn}
Of course if $\cU$ is a small open ball then both local systems~(\ref{sistloc}) are trivial.
Let $\Phi$ and $\Psi$ be markings of $\HH(V)|_{\cU}$ and $\HH(V^{\vee})|_{\delta(\cU)}$ respectively.
Then we have holomorphic local period maps
\begin{equation}\label{eccoper}
\begin{matrix}
\cU & \overset{\cP_{\Psi}}{\lra} & \cD_2 \\
A & \mapsto & \Phi_{\CC}(F^2_A)
\end{matrix}
\quad\quad
\begin{matrix}
\delta(\cU) & \overset{\cP_{\Phi}}{\lra} & \cD_2 \\
B & \mapsto & \Phi_{\CC}(F^2_{B}).
\end{matrix}
\end{equation}
The restrictions of $\cP_{\Psi}$ and $\cP_{\Phi}$ to $\cU\setminus\Delta^{\infty}_{*}(V)$ and $\delta(\cU)\setminus\Delta^{\infty}_{*}(V^{\vee})$ respectively
 are local period maps for  the families of double EPW-sextics parametrized by $\mathbb{LG}(\wedge^3 V)^0$ and  $\mathbb{LG}(\wedge^3 V^{\vee})^0$ respectively. We will be interested in comparing $\cP_{\Psi}(A)$ and $\cP_{\Phi}(A^{\bot})$ for $A\in \Delta^0_{*}(V)$ - thus $A^{\bot}\notin \mathbb{LG}(\wedge^3 V^{\vee})^0$.
Let $t\in T_{*}$ be such that $M_t\cong X_A$. By~(\ref{forzaitalia}) the marking $\Phi$ defines a marking $\Phi_{A^{\bot}}$ of $(S_t^{[2]},f_t^{*}\cO_{Y_B}(1))$ 
and 
\begin{equation}\label{valstrano}
\cP_{\Phi}(A^{\bot})=\cP_{\Phi_{A^{\bot}}}(S_t^{[2]},f_t^{*}\cO_{Y_B}(1)).
\end{equation}
 \section{Proof of Theorem~(\ref{princrisul})}
 \setcounter{equation}{0}
If $A\in \mathbb{LG}(\wedge^3 V)^{00}$ we have smooth double covers
$f_A\colon X_A\to Y_A$ and $f_{A^{\bot}}\colon X_{A^{\bot}}\to
Y_{A^{\bot}}$:  we will show that  $X_A$ and $X_{A^{\bot}}$ are
\lq\lq isogenous\rq\rq. Given a small open $\cU\subset\mathbb{LG}(\wedge^3 V)^{\sharp}$ we may consider  markings $\Psi$ and $\Phi$ of $\HH(V)|_{\cU}$ and $\HH(V^{\vee})|_{\delta(\cU)}$ respectively and the associated local period maps $\cP_{\Psi}$ and $\cP_{\Phi}$. 
We will show that
locally near $\Delta^0_{*}(V)$ we may choose $\Psi$ and $\Phi$ so that $\cP_{\Phi}\circ\delta$ is either
$r\circ\cP$ or the  composition of $r\circ\cP$ with a certain specific
reflection. In the final subsection we will rule out the latter case
 by considering the monodromy action; by analytic continuation this will prove Theorem~(\ref{princrisul}).
\subsection{Isogeny between $X_A$ and $X_{A^{\bot}}$}
\setcounter{equation}{0}
Let $A\in \mathbb{LG}(\wedge^3 V)^{00}$ and set $L_A:=f_A^{*}\cO_{Y_A}(1)$,  $L_{A^{\bot}}:=f_{A^{\bot}}^{*}\cO_{Y_{A^{\bot}}}(1)$. Let
\begin{align}
H^2(X_A)_{prim}:= & c_1(L_A)^{\bot}\subset H^2(X_A),\\
H^2(X_{A^{\bot}})_{prim}:= & c_1(L_{A^{\bot}})^{\bot}\subset H^2(X_{A^{\bot}}).
\end{align}
\begin{prp}\label{isogenia}
Let  $A\in \mathbb{LG}(\wedge^3 V)^{00}$. There exists an isomorphism of rational Hodge structures
\begin{equation}\label{eccogi}
g_A\colon H^2(X_A)_{prim}\overset{\sim}{\lra} H^2(X_{A^{\bot}})_{prim}
\end{equation}
well-defined up to $\pm 1$ and  such that for $\gamma\in
H^2(X_A)_{prim}$
\begin{equation}\label{interseco}
(\gamma,\gamma)_{X_A}= (g_A(\gamma),g_A(\gamma))_{X_{A^{\bot}}}.
\end{equation}
\end{prp}
\begin{proof}
We recall that $Y_{A^{\bot}}= Y_A^{\vee}$, see Corollary~(\ref{nontutto}).
Let $\Gamma_A\subset Y_A\times Y_{A^{\bot}}$ be the closure of the Gauss maps:
\begin{equation}
\Gamma_A:=\ov{\{(p,T_p Y_A)|\ p\in Y_A^{sm}\}}=
\ov{\{(T_q Y_{A^{\bot}},q)|\ q\in Y_{A^{\bot}}^{sm}\}}.
\end{equation}
Since $A\in\mathbb{LG}(\wedge^3 V)^0$
 the germ  of $Y_A$ at each of its
singular points is isomorphic to the product of $(\CC^2,0)$ and an $A_1$-singularity (see Proposition~(2.8) of~\cite{og2}) and similarly for $Y_{A^{\bot}}$. Thus the projection $\Gamma_A\to Y_A$ is identified with the blow-up $\wt{Y}_A\to Y_A$ of $sing(Y_A)$. Similarly the projection $\Gamma_A\to Y_{A^{\bot}}$ is identified with the blow-up $\wt{Y}_{A^{\bot}}\to Y_{A^{\bot}}$ of $sing(Y_{A^{\bot}})$.
 Thus $\Gamma_A$ defines an isomorphism $\wt{Y}_A\overset{\sim}{\lra} \wt{Y}_{A^{\bot}}$ and hence it gives an isomorphism of integral Hodge structures
\begin{equation}\label{accaquattro}
H^4(\wt{Y}_A)\overset{\sim}{\lra} H^4(\wt{Y}_{A^{\bot}}).
\end{equation}
The cohomology groups $H^4(\wt{Y}_A)$ and $H^4(X_A)$ are related as follows. Let $\phi_A\colon X_A\to X_A$ be the involution covering $f_A\colon X_A\to Y_A$. Let $F_A\subset X_A$ be the fixed locus of $\phi_A$: this is a Lagrangian smooth surface in $X_A$ because $\phi_A$ is anti-symplectic. Let $\wt{X}_A\to X_A$ be the blow-up of $F_A$; we have an isomorphism of integral Hodge structures
\begin{equation}
H^4(\wt{X}_A)\cong H^4(X_A)\oplus H^2(F_A)(-1).
\end{equation}
The involution $\phi_A$ lifts to an involution $\wt{\phi}_A\colon \wt{X}_A\to\wt{X}_A$ and $\wt{Y}_A\cong \wt{X}_A/\langle\wt{\phi}_A\rangle$. Thus we have an isomorphism of rational Hodge structures
\begin{equation}
H^4(\wt{Y}_A)\cong H^4(\wt{X}_{A})^{\langle\wt{\phi}_A\rangle}
\cong H^4(X_A)^{\langle\phi_A\rangle}\oplus H^2(F_A)(-1).
\end{equation}
Since $X_A$ is a deformation of $(K3)^{[2]}$ we have an isomorphism of rational Hodge Structures $Sym^2 H^2(X_A)\overset{\sim}{\lra} H^4(X_A)$ defined by cup-product. The action of $\phi_A$ on $H^2(X_A)$ has $(+1)$-eigenspace generated by $c_1(L_A)$ and $(-1)$-eigenspace equal to $H^2(X_A)_{prim}$ thus we get an isomorphism of rational Hodge structures
\begin{equation}
H^4(X_{A})^{\langle\phi_A\rangle}=
\CC c_1(L_A)^2\oplus Sym^2 H^2(X_A)_{prim}.
\end{equation}
The right-hand side of the above equality contains a rational $(2,2)$ class $q^{\vee}_A$ defined by \lq\lq inverting\rq\rq the Beauville-Bogomolov bilinear form (see Section~(3) of~\cite{og3}); let $W_A:=H^2(X_A)_{prim}\cap (q^{\vee}_A)^{\bot}$. One has a decomposition of rational Hodge Structures (Claim~(3.1) of~\cite{og3})
\begin{equation}
\CC c_1(L_A)^2\oplus \CC q^{\vee}_A\oplus
W_A.
\end{equation}
Of course we have analogous notions for $X_{A^{\bot}}$ and
hence~(\ref{accaquattro}) defines an isomorphism of rational H.S.'s
\begin{multline}\label{isolungo}
\CC c_1(L_A)^2\oplus \CC q^{\vee}_A\oplus W_A \oplus H^2(F_A)(-1)\cong \\
\CC c_1(L_{A^{\bot}})^2\oplus \CC q^{\vee}_{A^{\bot}}\oplus
W_{A^{\bot}}\oplus H^2(F_{A^{\bot}})(-1).
\end{multline}
If $A$ is very general (outside a countable union of proper analytic
subsets of $\mathbb{LG}(\wedge^3 V)^0
\cap\delta^{-1}\mathbb{LG}(\wedge^3 V^{\vee})^0$) then $W_A$ and
$W_{A^{\bot}}$ are both indecomposable rational Hodge Structures - see
Section~(3) of~\cite{og3}: since they contain the non-zero components
$H^{4,0}(X_A)$ and $H^{4,0}(X_{A^{\bot}})$ the above isomorphism
defines an isomorphism of rational H.S.'s
\begin{equation}
h_A\colon W_A\overset{\sim}{\lra} W_{A^{\bot}}.
\end{equation}
Let $Q_A\subset H^2(X_A)_{prim}$ and $Q_{A^{\bot}}\subset
H^2(X_{A^{\bot}})_{prim}$ be the cones of isotropic classes (with
respect to the Beauville-Bogomolov bilinear form). Let
 \begin{equation}
\begin{matrix}
H^2(X_A)_{prim} & \overset{\nu_A}{\lra} &
Sym^2 H^2(X_A)_{prim} \\
\alpha & \mapsto & \alpha^2
\end{matrix}
\end{equation}
be the (affine) Veronese map; we define similarly $\nu_{A^{\bot}}$.
Then
\begin{equation}
W_A=\text{span} \nu_A(Q_A),\quad
W_{A^{\bot}}=\text{span} \nu_{A^{\bot}}(Q_{A^{\bot}}).
\end{equation}
We claim that
\begin{equation}\label{nufissaqu}
h_A(\nu_A(Q_A))=\nu_{A^{\bot}}(Q_{A^{\bot}}).
\end{equation}
In fact let $\cU\subset \mathbb{LG}(\wedge^3 V)^0
\cap\delta^{-1}\mathbb{LG}(\wedge^3 V^{\vee})^0$ be an open ball
containing $A$. The Gauss-Manin connection gives  identifications
\begin{equation}\label{piattone}
H^2(X_B)\cong H^2(X_A),\quad
 H^2(X_{B^{\bot}})\cong H^2(X_{A^{\bot}})
\end{equation}
for every $B\in \cU$. Gauss-Manin gives also identifications $W_B\cong W_A$ and $W_{B^{\bot}}\cong W_{A^{\bot}}\cong$ sending $\nu_B(Q_B)$ to $\nu_A(Q_A)$ and $\nu_{B^{\bot}}(Q_{B^{\bot}})$ to $\nu_{A^{\bot}}(Q_{A^{\bot}})$ respectively.
The isomorphism $h_B$ is flat for the Gauss-Manin connection hence it  is identified with a (constant) map $h\colon W_A\to W_{A^{\bot}}$.
For $B\in \cU$ let $\sigma_B,\sigma_{B^{\bot}}$ be  symplectic forms on $X_B$ and $X_{B^{\bot}}$ respectively. Since $h_B$ is an isomorphism of H.S.'s we have
\begin{equation}\label{uffa}
h[\sigma_B^2]=[\sigma_{B^{\bot}}^2].
\end{equation}
Now $\sigma_B\in Q_A$ and $\sigma_{B^{\bot}}\in Q_{A^{\bot}}$ - here we
make the identifications~(\ref{piattone}) - and as $B$ varies in $\cU$
both $[\sigma_B]$ and $[\sigma_{B^{\bot}}]$ fill out  non-empty open
(in the classical topology) subsets of $\PP(Q_A)$ and
$\PP(Q_{A^{\bot}})$ respectively. Since $\PP(Q_A)$ and
$\PP(Q_{A^{\bot}})$ are  non singular quadrics any non-empty open
subset is Zariski-dense and hence~(\ref{uffa})
proves~(\ref{nufissaqu}). It follows from~(\ref{nufissaqu}) that there
exists a linear map~(\ref{eccogi}) well-defined up to $\pm 1$ such that
for $\alpha\in Q_A$ we have $h_A(\alpha^2)=g_A(\alpha)^2$.
By~(\ref{nufissaqu}) we have
\begin{equation}\label{montalbano}
\text{$g_A(\alpha)\in Q_{A^{\bot}}$ if and only if $\alpha\in Q_A$.}
\end{equation}
The rationality of $g_A$ follows from the fact that $h_A$ is defined
over $\QQ$, and $\nu_A,\nu_{A^{\bot}}$ give bijective maps between
$\PP(Q_A)(\QQ)$, $\nu_A(\PP(Q_A))(\QQ)$ and $\PP(Q_{A^{\bot}})(\QQ)$,
$\nu_{A^{\bot}}(\PP(Q_{A^{\bot}}))(\QQ)$ respectively. Finally let's
prove that Equation~(\ref{interseco}) holds. First we show that
\begin{equation}\label{crociato}
(\alpha,\beta)_{X_A}= (g_A(\alpha),g_A(\beta))_{X_{A^{\bot}}}, \quad
\alpha,\beta\in Q_A.
\end{equation}
Since~(\ref{isolungo}) respects the intersection forms we have
\begin{equation}\label{primaeq}
\int_{X_A}\alpha^2\wedge\beta^2=
\int_{X_{A^{\bot}}}h_A(\alpha^2)\wedge h_A(\beta^2)=
\int_{X_{A^{\bot}}}g_A(\alpha)^2\wedge g_A(\beta)^2.
\end{equation}
On the other hand $\alpha,\beta$ are isotropic and by~(\ref{nufissaqu})
$g_A(\alpha),g_A(\beta)$ are isotropic as well; thus
\begin{align}
\int_{X_A}\alpha^2\wedge\beta^2= & 2(\alpha,\beta)_{X_A}^2
\label{secondaeq}\\
\int_{X_{A^{\bot}}}g_A(\alpha)^2\wedge g_A(\beta)^2= &
 2 (g_A(\alpha),g_A(\beta))_{X_{A^{\bot}}}^2
\label{terzaeq}\\
\end{align}
(See Section~(2) of~\cite{og3}.) Equations~(\ref{primaeq}),
(\ref{secondaeq}) and (\ref{terzaeq}) prove that
either~(\ref{crociato}) holds or else
\begin{equation}\label{crociatomeno}
(\alpha,\beta)_{X_A}= -(g_A(\alpha),g_A(\beta))_{X_{A^{\bot}}}, \quad
\alpha,\beta\in Q_A.
\end{equation}
Assume that~(\ref{crociatomeno}) holds. Let $\sigma_A\in H^{2,0}(X_A)$
be the class of a symplectic forms on $X_A$; since $g_A$ is an
isomorphism of Hodge structures $g_A(\sigma_A)$ is represenyted by a
symplectic form on $X_{A^{\bot}}$. Then
\begin{equation}
(\sigma_A,\ov{\sigma_A})_{X_A}= -(g_A(\sigma_A),\ov{g_A(\sigma_A)})
_{X_{A^{\bot}}}.
\end{equation}
This is absurd because the Beaville-Bogomolov $(,)_X$ of an irreducible
symplectic manifold $X$ has the property that
$(\sigma,\ov{\sigma})_X>0$ for every symplectic form $\sigma$ on $X$.
This finishes the proof of Equation~(\ref{crociato}). Now let
$\gamma\in H^2(X_A)_{prim}$, then $\gamma=\alpha+\beta$ for certain
$\alpha,\beta\in Q_A$; by Equations~(\ref{montalbano})
and~(\ref{crociato} we get
\begin{equation}
(\gamma,\gamma)_{X_A}=
 2(\alpha,\beta)_{X_A}=2(g_A(\alpha),g_A(\beta))_{X_{A^{\bot}}}=
  (g_A(\gamma),g_A(\gamma))_{X_{A^{\bot}}}.
\end{equation}
This proves~(\ref{interseco}).
\end{proof}
Proposition~(\ref{isogenia}) implies that locally  there exists  $g\in O(\Lambda\otimes\QQ)$ which relates the periods of
$X_{A^{\bot}}$ to those of $X_A$ and furthermore $g$ normalizes the subgoup of monodromy operators in $O(\Lambda)$ . We  introduce some notation to
formalize this observation. 
Let $\cU\subset \mathbb{LG}(\wedge^3 V)^{\sharp}$ be a small open ball and let $\Psi,\Phi$ be markings of $\HH(V)|_{\cU}$ and $\HH(V^{\vee})|_{\delta(\cU)}$ respectively - see Definition~(\ref{marchio}).
Let $\ov{A}\in\cU$ be a reference point: the monodromy representation of
$\pi_1(\mathbb{LG}(\wedge^3 V)^{\sharp},\ov{A})$ on $\HH(V)$
determines via $\Psi_{\ov{A}}$ a  monodromy representation
$\pi_1(\mathbb{LG}(\wedge^3 V)^{\sharp},\ov{A})\to Stab(u)$ where $Stab(u)<
O(\wt{\Lambda})$ is the subgroup fixing the element $u$ given by~(\ref{scelgou}). An element of the image of the monodromy
representation is a {\it $\Psi_{\ov{A}}$-monodromy operator}. Similarly
we have a monodromy representation of $\pi_1(\mathbb{LG}(\wedge^3
V^{\vee}),\ov{A})$ on $\HH(V^{\vee})$;
 this determines via $\Phi_{\ov{A}^{\bot}}$ a  monodromy
representation $\pi_1(\mathbb{LG}(\wedge^3 V^{\vee})^{\sharp},\ov{A}^{\bot})\to
Stab(u)$. An element of the image of this second monodromy
representation is a {\it $\Phi_{\ov{A}^{\bot}}$-monodromy operator}. Here and
in the rest of the paper we will adopt the following conventions.
First we view both $O(\Lambda)$ and $O(\Lambda\otimes\QQ)$ as subgroups
of $O(\Lambda\otimes\CC)$. Secondly  if $\gamma\in
O(\Lambda\otimes\CC)$ we denote by $\wt{\gamma}\in
O(\wt{\Lambda}\otimes\CC)$ the isometry which fixes $u$ and equals
$\gamma$ on $\Lambda$. As a rule letters decorated by a tilde  denote
elements of $O(\wt{\Lambda}\otimes\CC)$, letters with no tilde denote
elements of $O(\Lambda\otimes\CC)$.
\begin{crl}\label{ortogonale}
Let $\cU\subset \mathbb{LG}(\wedge^3 V)^{\sharp}$ be a small open ball and let $\Psi,\Phi$ be markings of $\HH(V)|_{\cU}$ and $\HH(V^{\vee})|_{\delta(\cU)}$ respectively. Let $\cP_{\Psi},\cP_{\Phi}$ be the local period maps~(\ref{eccoper}). There exists $g\in
O(\Lambda\otimes\QQ)$ well-determined up to $\pm 1$ such that
\begin{equation}\label{operocongi}
\cP_{\Phi}(A^{\bot})=g\circ \cP_{\Psi}(A)
\end{equation}
for all $A\in\cU$. Let $\wt{\gamma}\in Stab(u)$ be a
$\Psi_{\ov{A}}$-monodromy operator; then
$\wt{g}\circ\wt{\gamma}\circ\wt{g}^{-1}$ is a $\Phi_{\ov{A}^{\bot}}$-monodromy
operator, in particular $g\circ\gamma\circ g^{-1}\in O(\Lambda)$.
\end{crl}
\begin{proof}
Equation~(\ref{operocongi}) holds on  $\cU\cap\mathbb{LG}(\wedge^3 V)^{00}$ by
Proposition~(\ref{isogenia}) and flatness of $g_A$; by
continuity Equation~(\ref{operocongi}) holds on all of $\cU$. The statement about monodromy operators holds by flatness of $g_A$.
\end{proof}
\subsection{Restriction to $\Delta^0_{*}(V)$ of the
local period maps}
\setcounter{equation}{0}
Our next task is to analyze the restriction to $\Delta^0_{*}(V)$ of the
local period maps $\cP_{\Psi}$ and $\cP_{\Phi}$.
\begin{prp}\label{buonimarchi}
Let $\cU\subset \mathbb{LG}(\wedge^3 V)^{\sharp}$ be a small open ball.  There is a choice of markings
$\Psi$ and $\Phi$  of $\HH(V)|_{\cU}$ and $\HH(V^{\vee})|_{\delta(\cU)}$ respectively    such that
\begin{equation}\label{divisore}
\cP_{\Psi}(\Delta^0_{*}(V)\cap\cU)=(e_1+2e_2)^{\bot}\cap\cP_{\Psi}(\cU)
\end{equation}
where $e_1,e_2\in\Lambda$ are as in Section~(\ref{prologo}). Furthermore
\begin{equation}\label{bellaformula}
\cP_{\Phi}(A^{\bot})=r\circ\cP_{\Psi}(A),\quad A\in \Delta^0_{*}(V)
\end{equation}
where $r$ is the involution defined by~(\ref{alfaexpl}).
\end{prp}
\begin{proof}
First we embed the lattice $\wt{\Lambda}$ in a unimodular lattice as follows. Let $\wh{\Lambda}:=U^4\wh{\oplus}(-E_8)^2$. Let $U_1<\wh{\Lambda}$ be one of the hyperbolic lattices, let $z\in U_1$ be a vector of square $2$ and $e_2$ be a generator of $z^{\bot}\cap U_1$. Then  we have an isomorphism
\begin{equation}\label{zetaperp}
z^{\bot}\cong \wt{\Lambda}
\end{equation}
and we can choose it so that it matches the present $e_2$ with the vector $e_2$ appearing in~(\ref{lambdazero}): we fix such an isomorphism once for all. Let $u,e_1\in\wt{\Lambda}=z^{\bot}$ be as in
Section~(\ref{prologo}); then $\{(u\pm e_1)/2,(z\pm e_2)/2\}\subset \wh{\Lambda}$. Furthermore the sublattices $\langle (u+ e_1)/2,(u- e_1)/2\rangle$ and $\langle(z+ e_2)/2,(z- e_2)/2 \rangle$ are orthogonal hyperbolic planes. Thus we have an orthogonal decomposition
\begin{equation}\label{inzaghi}
\wh{\Lambda}=\langle (u+ e_1)/2,(u- e_1)/2\rangle\wh{\oplus}
\langle(z+ e_2)/2,(z- e_2)/2 \rangle\wh{\oplus} U^2
\wh{\oplus} (-E_8)^2.
\end{equation}
Now we pass to the geometry.
Let $T_{*}$ be as in Subsection~(\ref{dualesestica}). Let
$A\in\Delta^0_{*}(V)\cap\cU$; by definition there exists $t\in T_{*}$ such
that
\begin{equation}\label{eccovu}
M_t\cong X_{A}.
\end{equation}
Since $\cU$ is a small open ball there exists a small open ball
$\cV\subset T_{*}$ such that if $t\in\cV$ then~(\ref{eccovu})
holds for some $A\in\Delta^0_{*}(V)\cap\cU$ and conversely if $A\in\Delta^0_{*}(V)\cap\cU$ then there exists $t\in\cV$ such that~(\ref{eccovu})
holds. Let $\kappa\colon\cS\to T_{*}$ be the tautological family of $K3$
surfaces parametrized by $T_{*}$. Let $t\in T_{*}$. Let $S_t=\kappa^{-1}(t)$, $D_t$ etc.~be as in Subsection~(\ref{deltazero}) and $A_t\in\Delta^0_{*}(V)\cap\cU$ such that $M_t\cong X_{A_t}$.
Let $f_t\colon M_t\to Y_{A_t}$ be the double cover.
Let $v_t,w_t\in  H^{*}(S_t)$ be given by
\begin{equation}
v_t:= 2+c_1(D_t)+2\eta_t,\quad w_t:=1-\eta_t
\end{equation}
where $\eta_t\in H^4(S_t;\ZZ)$ is the orientation class. The Mukai map
\begin{equation}\label{vuperp}
\theta_{v_t}\colon v_t^{\bot}\to H^2(M_t)
\end{equation}
is an isometry of Hodge structures - see Subsection~(\ref{deltazero}).
Furthermore one has
\begin{equation}\label{taddei}
c_1(f_t^{*}L_{A_t})=\theta_{v_t}(\eta_t-1)=\theta_{v_t}(-w_t^{\vee}).
\end{equation}
(See the line preceding~(\ref{mukform}) for the definition of
$w_t^{\vee}$.)  The local system $R^2\kappa_{*}\ZZ|_{\cV}$ is trivial because $\cV$ is a small open ball. Thus
there exist  sections $\alpha,\beta\in\Gamma(R^2\kappa_{*}\ZZ|_{\cV})$
such that $c_1(D_t)=\alpha_t+5\beta_t$ for all $t\in\cV$. We define a
trivialization
\begin{equation}\label{noia}
R\kappa_{*}\ZZ|_{\cV}= (R^0\kappa_{*}\ZZ\oplus R^2\kappa_{*}\ZZ\oplus
R^4\kappa_{*}\ZZ)|_{\cV} \overset{\Upsilon}{\lra} \cV\times\wh{\Lambda}
\end{equation}
as follows.
For $t\in\cV$ let
\begin{equation}
\begin{array}{rcclclclcl}
\Upsilon_t(1)  & :=  & - & u/2  &+ & e_1/2 & - & z  & + & e_2 \\
\Upsilon_t(\eta_t)  & :=&   & u/2  &+& e_1/2 & - & z & + & e_2 \\
\Upsilon_t(\alpha_t)  &  :=&  &  & - & 2e_1 &+ & 5z/2 & - & 3e_2/2   \\
\Upsilon_t(\beta_t)  & :=& & & &   &  & z/2 & - & e_2/2  \\
\end{array}
\end{equation}
and let
\begin{equation}
\Upsilon_t|_{\{1,\eta_t,\alpha_t,\beta_t\}^{\bot}}
\colon \{1,\eta_t,\alpha_t,\beta_t\}^{\bot}
\overset{\sim}{\lra}
\{(u\pm e_1)/2,(z\pm e_2)/2\}^{\bot}
\end{equation}
be an arbitrary isometry - notice that
$\{1,\eta_t,\alpha_t,\beta_t\}^{\bot}\cong U^2\wh{\oplus}(-E_8)^2$ is
isometric to $\{(u\pm e_1)/2,(z\pm e_2)/2\}^{\bot}$ by~(\ref{inzaghi}).
A straightforward computation shows that $\Upsilon_t$ is  an isometry.
 The trivialization~(\ref{noia}) is defined to have value $\Upsilon_t$
 at $t\in\cV$. Now notice that $\Upsilon_t(v_t)=z$ and hence we have an
isometry $\Upsilon_t\circ\theta_{v_t}^{-1}\colon
H^2(M_t)\overset{\sim}{\lra}z^{\bot}=\wt{\Lambda}$. Since
 $\Upsilon_t(-w_t^{\vee})=u$ Equation~(\ref{taddei}) gives that
\begin{equation}
\Upsilon_t\circ\theta_{v_t}^{-1}(c_1(f_t^{*}L_{A_t}))=u.
\end{equation}
Hence $\Upsilon_t\circ\theta_{v_t}^{-1}$ defines a marking of
$(M_t,f_t^{*}L_{A_t})$ for every $t\in\cV$;  since
$\HH(V)|_{\cU}$ is trivial there exists a marking $\Psi\colon
\HH(V)|_{\cU}\lra\cU\times\wt{\Lambda}$ such that
\begin{equation}
\Psi_{A_t}=\Upsilon_t\circ\theta_{v_t}^{-1},\quad t\in\cV. 
\end{equation}
Equation~(\ref{divisore}) follows from Equation~(\ref{imagodelta}) and the equality
\begin{equation}\label{epiuduee}
\Psi_{A_t}(\theta_{v_t}(5+2 c_1(D_t)+5\eta_t))
 =\Upsilon_t(5+2 c_1(D_t)+5\eta_t)=e_1+2e_2.
\end{equation}
 Next we define a
marking $\Phi$ for  $\HH(V^{\vee})|_{\delta(\cU)}$. 
By~(\ref{forzaitalia}) this will be equivalent to a marking of $S_t^{[2]}$, hence we first recall the
description of $H^2(S_t^{[2]})$. We notice that $w_t$ is the Mukai
vector (see~(\ref{vettoremukai})) of any ideal sheaf $I_Z$ where
$[Z]\in S_t^{[2]}$; Mukai's map
\begin{equation}\label{mukaiwt}
\theta_{w_t}\colon w_t^{\bot}\to H^2(S_t^{[2]})
\end{equation}
is an isomorphism of polarized Hodge structures and an isometry. Let
$g_t\colon S_t^{[2]}\to Y_{A_t^{\bot}}$ be the map defined
by~(\ref{mappaeffe}); by Subsection~(5.3) of~\cite{og1}
\begin{equation}\label{aquilani}
c_1(g_t^{*}\cO_{Y_{A_t^{\bot}}}(1))=\theta_{w_t}( -2+c_1(D_t)-2\eta_t)
=\theta_{w_t}(-v_t^{\vee}).
\end{equation}
We define a trivialization
\begin{equation}\label{noiadue}
R\kappa_{*}\ZZ|_{\cV}=
(R^0\kappa_{*}\ZZ\oplus R^2\kappa_{*}\ZZ\oplus R^4\kappa_{*}\ZZ)|_{\cV}\overset{\Theta}{\lra} \cV\times\wh{\Lambda}
\end{equation}
as follows.
For $t\in\cV$ let
\begin{equation}
\begin{array}{rclclclcl}
\Theta_t(1)  & :=  &   u  &- & e_1 & + & z/2  & - & e_2/2 \\
\Theta_t(\eta_t)  & :=&    u  & - & e_1 & - & z/2 & - & e_2/2 \\
\Theta_t(\alpha_t)  &  :=&    5u/2  & - & 3e_1/2 & &  & - & 2e_2 \\
\Theta_t(\beta_t)  & :=&  u/2 & - & e_1/2  &  &  &  &   \\
\end{array}
\end{equation}
and let the restriction of $\Theta_t$ to $\{1,\eta_t,\alpha_t,\beta_t\}^{\bot}$ be equal to the restriction of $\Upsilon_t$.
A straightforward computation shows that $\Theta_t$ is an isometry.
The trivialization~(\ref{noiadue}) is defined to have value $\Theta_t$ at $t\in\cV$. Now notice that $\Theta_t(w_t)=z$
and hence we have an isometry $\Theta_t\circ\theta_{w_t}^{-1}\colon H^2(S^{[2]}_t)\overset{\sim}{\lra}z^{\bot}=\wt{\Lambda}$. Since
 $\Theta_t(-v_t^{\vee})=u$ Equation~(\ref{aquilani}) gives that
\begin{equation}
\Theta_t\circ\theta_{v_t}^{-1}(c_1(g_t^{*}\cO_{Y_{A^{\bot}_t}}))=u.
\end{equation}
Hence $\Theta_t\circ\theta_{w_t}^{-1}$ defines a marking of
$(S_t^{[2]},g_t^{*}\cO_{Y_{A^{\bot}_t}}(1))$ for every $t\in\cV$;  by~(\ref{forzaitalia}) and triviality of
$\HH(V^{\vee})|_{\delta(\cU)}$  there exists a marking $\Phi\colon
\HH(V^{\vee})|_{\delta(\cU)}\lra\cU\times\wt{\Lambda}$ such that
\begin{equation}
\Phi_{A_t}=\Theta_t\circ\theta_{w_t}^{-1},\quad t\in\cV. 
\end{equation}
 Now let's
prove~(\ref{bellaformula}).  Equation~(\ref{bellaformula})
is equivalent to
\begin{equation}
\cP_{\Theta_t\circ\theta_{w_t}^{-1}}
(S_t^{[2]},g_t^{*}\cO_{Y_{A^{\bot}_t}}(1))=
r\circ\cP_{\Upsilon_t\circ\theta_{v_t}^{-1}}(M_t,f_t^{*}L_A),
\quad t\in\cV
\end{equation}
by Equation~(\ref{valstrano}).
Since $\theta_{v_t}$ and $\theta_{w_t}$ are isomorphism of Hodge structures the above equation may be rewitten as
\begin{equation}\label{storace}
\Theta_t(H^{2,0}(S_t))=
r\circ\Upsilon_t(H^{2,0}(S_t)),\quad t\in\cV.
\end{equation}
For $t\in\cV$ let
\begin{equation}
\begin{matrix}
H^{*}(S_t;\QQ) & \overset{\Xi_t}{\lra} & H^{*}(S_t;\QQ) \\
\gamma & \mapsto & -\gamma^{\vee}+
\frac{1}{2}(\gamma^{\vee},v_t+w_t)(v_t+w_t)
\end{matrix}
\end{equation}
i.e.~the composition of the isometry $\gamma\mapsto\gamma^{\vee}$ and the reflection which is $(+1)$ on $\QQ (v_t+w_t)$ and $(-1)$ on $(v_t+w_t)^{\bot}$; thus $\Xi_t$ is a rational (not integral !) Hodge isometry. A straightforward computation shows that $\Theta_t=r\circ\Upsilon_t\circ\Xi_t$. Thus~(\ref{storace}) follows at once from $\Xi_t (H^{2,0}(S))=H^{2,0}(S)$.
\end{proof}
If $\gamma\in\Lambda\otimes\QQ$ is non-isotropic we let
\begin{equation}
\begin{matrix}
\Lambda\otimes \QQ & \overset{r_{\gamma}}{\lra} & \Lambda\otimes\QQ \\
x & \mapsto & -x+\frac{2}{(\gamma,\gamma)}(x,\gamma)\gamma
\end{matrix}
\end{equation}
be the reflection with $(+1)$-eigenspace $\QQ\gamma$ and
$(-1)$-eigenspace $\gamma^{\bot}$.
Let
\begin{equation}
\zeta:= e_1+2e_2.
\end{equation}
\begin{crl}\label{autaut}
Keep notation and assumptions of Proposition~(\ref{buonimarchi}). Let $\Psi$, $\Phi$ be the markings of Proposition~(\ref{buonimarchi}). Then:
\begin{itemize}
\item[(1)]
$\cP_{\Phi}(A^{\bot})=r\circ\cP_{\Psi}(A)$ for all $A\in\cU$ or
\item[(2)]
$\cP_{\Phi}(A^{\bot})=r\circ r_{\zeta}\circ\cP_{\Psi}(A)$ for all
$A\in\cU$.
\end{itemize}
\end{crl}
\begin{proof}
By Corollary~(\ref{ortogonale}) there exists $g\in O(\Lambda\otimes\QQ)$ such that $\cP_{\Phi}(A^{\bot})=g\circ\cP_{\Psi}(A)$ for all
$A\in\cU$.
By Proposition~(\ref{buonimarchi}) $r^{-1}\circ g$ fixes the points of $\zeta^{\bot}\cap\cD_2$: it follows that $r^{-1}\circ g$ fixes $\zeta^{\bot}\subset\PP(\Lambda\otimes\CC)$. Since $r\in O(\Lambda)$ we have
$r^{-1}\circ g\in O(\Lambda\otimes\QQ)$: thus we get that
 that
\begin{equation}
r^{-1}\circ g|_{\zeta^{\bot}}=\pm Id_{\zeta^{\bot}}.
\end{equation}
Since $\zeta$ is non-isotropic (in fact $(\zeta,\zeta)=-10$) and $r^{-1}\circ g\in O(\Lambda\otimes\QQ)$
we get that $r^{-1}\circ g(\zeta)=\pm\zeta$. It follows that  $r^{-1}\circ g=\pm Id$ or $r^{-1}\circ g=\pm r_{\zeta}$.  Since $-Id$ acts trivially on $\cD_2\subset\PP(\Lambda)$ the corollary follows.
\end{proof}
\subsection{The proof}
\setcounter{equation}{0}
We will apply the monodromy statement of Corollary~(\ref{ortogonale}) in order to show that Item~(2) of Corollary~(\ref{autaut}) can not hold. We will use the following result.
\begin{clm}\label{conto}
Let
\begin{equation}
\xi=a_1e_1+a_2 e_2+\nu\in\Lambda
\end{equation}
be a $(-2)$-vector i.e.~$(\xi,\xi)=-2$ and assume that $r_{\zeta}\circ
r_{\xi}\circ r_{\zeta}\in O(\Lambda)$. Then
\begin{equation}\label{cinque}
(\xi,\zeta)\equiv 0 \pmod{5}.
\end{equation}
\end{clm}
\begin{proof}
A tedious
straightforward computation gives that
\begin{multline}\label{rifdieuno}
r_{\zeta}\circ r_{\xi}\circ r_{\zeta}(e_1)=
\frac{1}{25}
(18 a_1^2-48 a_1 a_2 + 32a_2^2-25)e_1-\\
-\frac{1}{25}(24 a_1^2-14 a_1 a_2 -24 a_2^2) e_2 -
\frac{2}{5}(3a_1-4a_2)\nu.
\end{multline}
Thus $18 a_1^2-48 a_1 a_2 + 32a_2^2\equiv 0\pmod{25}$. Since
\begin{equation}
18 a_1^2-48 a_1 a_2 + 32a_2^2= 2(3a_1-4a_2)^2
\end{equation}
we get that $3a_1-4a_2\equiv 0\pmod{5}$. This proves~(\ref{cinque}) because
\begin{equation}
 (\xi,\zeta)=-2a_1-4a_2\equiv 3a_1-4a_2\pmod{5}.
\end{equation}
\end{proof}
\begin{prp}\label{rifmon}
Keep notation and assumptions of Proposition~(\ref{buonimarchi}). Let
$\Psi$, $\Phi$ be the markings of Proposition~(\ref{buonimarchi}) and
$\ov{A}\in\cU\cap\Delta^0_{*}(V)$. There exists a $(-2)$-vector
$\xi\in\Lambda$ such that $-r_{\xi}$ is a $\Psi_{\ov{A}}$-monodromy
operator and
\begin{equation}\label{nondivisibile}
(\xi,\zeta)\not\equiv 0\pmod{5}.
\end{equation}
\end{prp}
We grant the above proposition for the moment being and we proceed to
prove Theorem~(\ref{princrisul}). Let notation and assumptions be as in
Proposition~(\ref{buonimarchi}) and $\Psi$, $\Phi$ be the markings of
Proposition~(\ref{buonimarchi}). Then either~(1) or~(2) of
Corollary~(\ref{autaut}) holds.  Suppose that~(2) holds; we will arrive
at a contradiction. Let $\xi$ be as in Proposition~(\ref{rifmon}): by
Corollary~(\ref{ortogonale}) we have $-r\circ r_{\zeta}\circ
r_{\xi}\circ r_{\zeta}\circ r\in O(\Lambda)$. Since $r\in O(\Lambda)$
we get that $r_{\zeta}\circ r_{\xi}\circ r_{\zeta}\in O(\Lambda)$: this
contradicts Claim~(\ref{conto}) because of~(\ref{nondivisibile}).
Thus~(1) of Corollary~(\ref{autaut}) holds. Let
$\cU^{00}:=\cU\cap\mathbb{LG}(\wedge^3 V)^{00}$; then $\cU^{00}$ is an
open (in the euclidean topology) non-empty subset of
$\mathbb{LG}(\wedge^3 V)^{00}$. Since~(1) of Corollary~(\ref{autaut})
holds we have
\begin{equation}
\cP\circ\delta|_{\cU^{00}}=\overline{r}\circ\cP|_{\cU^{00}}.
\end{equation}
Both $\cP\circ\delta$ and $\overline{r}\circ\cP$ are holomorphic maps
with domain the connected manifold $\mathbb{LG}(\wedge^3 V)^{00}$; by
analytic continuation we get that Theorem~(\ref{princrisul}) holds.
\vskip 2mm
 \noindent
 {\bf Proof of Proposition~(\ref{rifmon})}
 \hskip 2mm
  Let $F\subset\PP^3$ be a smooth quartic, thus $F$ is a $K3$
surface. We have a regular map
\begin{equation}\label{mappagi}
\begin{matrix}
F^{[2]} & \overset{g}{\lra} & {\mathbb Gr}(1,\PP^3)\subset\PP^5\\
[Z] & \mapsto & span(Z)
\end{matrix}
\end{equation}
 and $c_1(g^{*}\cO_{ {\mathbb Gr}(1,\PP^3)}(1))$ has square $2$ for
the Beauville-Bogomolov form. If $F$ does not contain lines the above
map is finite and hence $g^{*}\cO_{ {\mathbb Gr}(1,\PP^3)}(1)$ is an
ample line-bundle on $F^{[2]}$, if $F$ contains a line $R$ then
$g^{*}\cO_{ {\mathbb Gr}(1,\PP^3)}(1)$ is big and nef but it restricts
to the trivial line-bundle on the $\PP^2$ given by $R^{(2)}\subset
F^{[2]}$. Assume that $F$ does not contain lines; we proved in
Section~(6) of~\cite{og1} that
\begin{equation}\label{sondef}
\text{$(F^{[2]},g^{*}\cO_{ {\mathbb Gr}(1,\PP^3)}(1))$ is deformation
equivalent to $(M_t,f_t^{*}L_{A_t})$}
\end{equation}
where $t\in T$, i.e.~there exists a polarized family of irreducible
symplectic $4$-folds over a connected basis with one fiber isomorphic
to $(F^{[2]},g^{*}\cO_{ {\mathbb Gr}(1,\PP^3)}(1))$ and another fiber
isomorphic to $(M_t,f_t^{*}L_{A_t})$. Using this result we will show
that the monodromy operator on $F$ given by a suitable $(-2)$-class
orthogonal to $c_1(\cO_F(1))$ gives rise to a $\Psi(\ov{A})$-monodromy
operator for which~(\ref{nondivisibile}) holds. Before proving this we
must dive into the details of the proof of~(\ref{sondef}). Let
$F_0\subset\PP^3$ be a smooth quartic surface containing a line $R$ and
with Picard number $2$, i.e.~$Pic(F_0)=\ZZ [A_0]\oplus\ZZ[R]$ where
$A_0$ is the (hyper)plane class. The divisor $(2A_0-R)$ is very ample
and $c_1(2A_0-R)^2=10$; thus we have an embedding
\begin{equation}\label{effeimmersa}
F_0\hra |2A_0-R|^{\vee}\cong\PP^6
\end{equation}
as a linearly normal $K3$ surface of degree $10$. Let
\begin{equation}
w_0:=1+c_1(A_0)+\eta_0,\quad v_0:=2+c_1(2A_0-R)+2\eta_0
\end{equation}
where $\eta_0\in H^4(F_0;\ZZ)$ is the orientation class. Let $M_{w_0}$
be the moduli space of torsion-free sheaves $\cG$ on $F_0$ such that
$v(\cG)=w_0$; every such sheaf is equal to $I_{Z}\otimes\cO_{F_0}(A_0)$
for a unique $[Z]\in F_0^{[2]}$ and hence
\begin{equation}\label{stessi}
M_{w_0}\cong F_0^{[2]}.
\end{equation}
We let $L_{w_0}:=g_0^{*}\cO_{ {\mathbb Gr}(1,\PP^3)}(1)$. Let $M_{v_0}$
be the moduli space of $(2A_0-R)$-semistable sheaves $\cF$ on $F_0$
such that $v(\cF)=v_0$. The moduli space $M_{v_0}$ is smooth because
$(2A_0-R)$ is $v_0$-generic (see Section~(6) of~\cite{og1}) and hence
it is a deformation of $(K3)^{[2]}$. Mukai's map gives isometries of
Hodge structures
\begin{equation}
\theta_{w_0}\colon w_0^{\bot}\overset{\sim}{\lra}H^2(M_{w_0}),\quad
\theta_{v_0}\colon v_0^{\bot}\overset{\sim}{\lra}H^2(M_{v_0}).
\end{equation}
One has  (see p.1241 of~\cite{og1})
\begin{equation}\label{ellewu}
c_1(L_{w_0})=\theta_{w_0}(\eta_0-1).
\end{equation}
We let $L_{v_0}$ be the line-bundle on $M_{v_0}$ such that
\begin{equation}\label{ellevu}
c_1(L_{v_0})=\theta_{v_0}(\eta_0-1).
\end{equation}
In Lemma~(6.2) of~\cite{og1} we considered the birational map
$M_{w_0}\dashrightarrow M_{v_0}$ whose inverse
\begin{equation}\label{inversadi}
\varphi\colon M_{v_0}\dashrightarrow M_{w_0}
\end{equation}
is the Mukai reflection defined by the $(-2)$-vector
\begin{equation}\label{uzero}
u_0:=(1+c_1(A_0-R)+\eta_0)
\end{equation}
(notice that $-r_{u_0}(v_0)=w_0$), see~\cite{yoshi}. Since $M_{v_0}$
and $M_{w_0}$ are irreducible symplectic manifolds the birational map
$\varphi$ induces an isomorphism of lattices $\varphi^{*}\colon
H^2(M_{w_0})\cong H^2(M_{v_0})$. By Theorem~(2.9) of~\cite{yoshi} we
have
\begin{equation}\label{riflesso}
\varphi^{*}\theta_{w_0}(\alpha)= \theta_{v_0}(-r_{u_0}(v_0)),
\end{equation}
in particular by~(\ref{ellewu})-(\ref{ellevu}) we have
\begin{equation}\label{corrispondono}
\varphi^{*}L_{w_0}\cong L_{v_0}.
\end{equation}
The birational map $\varphi$ is the flop of
\begin{equation}
\Pi_{w_0}:=R^{(2)}\subset F_0^{[2]}=M_{w_0}.
\end{equation}
It follows that $M_{v_0}$ contains $\Pi_{v_0}\cong\Pi_{w_0}^{\vee}$ and
from~(\ref{corrispondono}) we get that $L_{v_0}$ is big, nef and its
restriction to $\Pi_{v_0}$ is trivial.  Let $\cX\to B_{v_0}$  be a
representative for the deformation space of $(M_{v_0},L_{v_0})$,
i.e.~deformations of $M_{v_0}$ that \lq\lq keep $c_1(L_{v_0})$ of type
$(1,1)$\rq\rq. Similarly let $\cX'\to B_{w_0}$ be a representative for
the deformation space of $(M_{w_0},L_{w_0})$. We let $0\in B_{v_0}$ and
$0\in B_{w_0}$ be the points corresponding to $(M_{v_0},L_{v_0})$ and
$(M_{w_0},L_{w_0})$ respectively. Thus for each $q\in B_{v_0}$ the
fiber $X_q$ of $\cX\to B_{v_0}$ over $q$ has a line-bundle $L_q$ which
is a deformation of $L_{v_0}$. Similarly for each $s\in B_{w_0}$ the
fiber $X'_s$ of $\cX'\to B_{w_0}$ over $s$ has a line-bundle $L'_s$
which is a deformation of $L_{w_0}$. We may and will assume that
$B_{v_0},B_{w_0}$ are contractible and hence Gauss-Manin gives
identifications
\begin{equation}\label{banale}
H^2(X_q)\cong H^2(M_{v_0}),\quad H^2(X'_s)\cong H^2(M_{w_0}),
 \quad q\in B_{v_0},\ s\in B_{w_0}
\end{equation}
which match $c_1(L_q)$ to $c_1(L_{v_0})$ and $c_1(L'_s)$ to
$c_1(L_{w_0})$. Let $B(\Pi_{v_0})\subset B_{v_0}$ be the locus
parametrizing deformations $X_q$ which contain a deformation of
$\Pi_{v_0}$, and similarly let $B(\Pi_{w_0})\subset B_{w_0}$ be the
locus parametrizing deformations of $X'_s$ which contain a deformation
of $\Pi_{w_0}$. By a Theorem of Voisin~\cite{voi} each of these loci is
smooth of codimension $1$. There is a natural isomorphism of germs
\begin{equation}
\mu\colon (B_{w_0},0)\overset{\sim}{\lra}(B_{v_0},0)
\end{equation}
such that $\mu(\Pi_{w_0})=\Pi_{v_0}$ and if $s\notin \Pi_{w_0}$ then
$(X_{\mu(s)},L_{\mu(s)})\cong (X'_{s},L'_{s})$. Let $s\in B_{w_0}$ be
such that $(X'_{s},L'_{s})\cong (F^{[2]},g^{*}\cO_{ {\mathbb
Gr}(1,\PP^3)}(1))$ where $F$ is a quartic containing no lines. Then
$s\notin \Pi_{w_0}$ and hence
\begin{equation}\label{isomchiave}
(X_{\mu(s)},L_{\mu(s)})\cong (X'_{s},L'_{s})
 \cong (F^{[2]},g^{*}\cO_{{\mathbb Gr}(1,\PP^3)}(1)).
\end{equation}
 On the other hand there
exists $q\in B_{v_0}$ such that $(X_q,c_1(L_q))\cong
(M_t,\theta_t(\eta_t-1))$ for $t\in T_{*}$ because the parameter space
for linearly normal $K3$ surfaces of degree $10$ (an open subset of the
relevant Hilbert scheme) is irreducible - notice that
$q\notin\Pi_{v_0}$ because $\theta_t(\eta_t-1)$ is ample. Thus there
exists $\ov{A}\in \Delta^0_{*}(V)$ such that
\begin{equation}\label{altrachiave}
(X_{\ov{A}},L_{\ov{A}})\cong(X_q,L_q).
\end{equation}
Let $\gamma\in H^2(F;\ZZ)$ be a $(-2)$-class orthogonal to
$c_1(\cO_F(1))$.  Then $\gamma$ determines a monodromy operator
 on $H^2(F^{[2]})$; by~(\ref{isomchiave}) and~(\ref{altrachiave})  this
 monodromy operator can be identified with a monodromy operator on
 $X_{\ov{A}}$ because polarized deformation spaces of  irreducible symplectic
 manifolds are smooth. Given the trivializations~(\ref{banale})  the
 monodromy operator in question is equal to
$-r_{\varphi^{*}\theta_{w_0}(\gamma_0)}$ - here $\gamma_0\in H^2(F_0)$
is the class corresponding to $\gamma$ (the second trivialization
of~(\ref{banale}) defines an identification $H^2(F)\cong H^2(F_0)$).
The corresponding $\Psi_{\ov{A}}$-monodromy operator is equal to
$-r_{\xi}$ where
\begin{equation}\label{eccoxi}
\xi:=\Psi_{\ov{A}}(\varphi^{*}\theta_{w_0}(\gamma_0))\in\Lambda
\end{equation}
\begin{clm}\label{ultimo}
Keep notation as above. Let $\gamma_0\in c_1(A_0)^{\bot}\subset
H^2(F_0;\ZZ)$ be a $(-2)$-class such that
\begin{equation}\label{noncinque}
\int_{F_0}\gamma_0\wedge c_1(R)\not\equiv 0\pmod{5}.
\end{equation}
Then
\begin{equation}\label{finalmente}
(\xi,\zeta)\not\equiv 0\pmod{5}.
\end{equation}
\end{clm}
\begin{proof}
By~(\ref{epiuduee}) we must check that
\begin{equation}\label{daverificare}
(\varphi^{*}\theta_{w_0}(\gamma_0),\theta_{v_0}(5+2
c_1(2A_0-R)+5\eta_0))
 \not\equiv 0\pmod{5}.
\end{equation}
By~(\ref{riflesso}) this is equivalent to
\begin{equation}
(-r_{u_0}(\gamma_0),5+2 c_1(2A_0-R)+5\eta_0)\not\equiv 0\pmod{5}.
\end{equation}
where $u_0$ is given by~(\ref{uzero}). A straightforward computation
gives that
\begin{equation}
(-r_{u_0}(\gamma_0),5+2 c_1(2A_0-R)+5\eta_0)=
-8\int_{F_0}\gamma_0\wedge c_1(R)
\end{equation}
and hence we get~(\ref{daverificare}).
\end{proof}
Certainly there exists a class $\gamma_0$ satisfying the hypotheses of
Claim~(\ref{ultimo}): the vector $\xi\in\Lambda$ given
by~(\ref{eccoxi}) satisfies the thesis of Proposition~(\ref{rifmon}).
\vskip 1cm
 \scriptsize{
Kieran G. O'Grady\\
Universit\`a di Roma ``La Sapienza",\\
Dipartimento di Matematica ``Guido Castelnuovo",\\
Piazzale Aldo Moro n.~5, 00185 Rome, Italy,\\
e-mail: {\tt ogrady@mat.uniroma1.it}. }

\end{document}